\documentclass[final]{elsarticle}

\usepackage{amssymb}
\usepackage{amsmath}
\usepackage{amsthm}
\usepackage{algorithm,algpseudocode}
\usepackage{hyperref}
\hypersetup{colorlinks=true}
\usepackage[dvipsnames]{xcolor}

\DeclareMathOperator*{\argmin}{arg\,min}

\newcommand{\id}{\mathrm d}

\newcommand{\vc}{\mathbf}

\renewcommand{\tilde}{\widetilde}
\newcommand{\pard}[2]{\frac{\partial #1}{\partial #2}}

\newcommand{\fm}{\varphi}

\newcommand{\err}{\mathcal E}
\newcommand{\bphi}{\pmb\phi}
\newcommand{\balpha}{\pmb\alpha}

\newtheorem{thm}{Theorem}

\newtheorem{lem}{Lemma}

\newtheorem{defn}{Definition}
\newtheorem{ass}{Assumption}

\journal{Journal of Computational Science}

\begin{document}

\begin{frontmatter}
\cortext[cor1]{Corresponding author.}

\title{Discrete Empirical Interpolation Method with Upper and Lower Bound Constraints}
\author[aff1]{Louisa B. Ebby}
\ead{lbebby@ncsu.edu}
\author[aff1]{Mohammad Farazmand\corref{cor1}}
\ead{farazmand@ncsu.edu}
        \address[aff1]{Department of Mathematics, North Carolina State University,
        	2311 Stinson Drive, 
        	Raleigh,
        	27695-8205, 
        	NC,
        	USA}

\begin{abstract}
Discrete Empirical Interpolation Method (DEIM) is a simple and effective method for reconstructing a function from its incomplete pointwise observations.
However, applying DEIM to functions with physically constrained ranges can produce reconstructions with values outside the prescribed physical bounds. 
Such physically constrained quantities occur routinely in applications, e.g., mass density whose range is nonnegative.
The DEIM reconstructions which violate these physical constraints are not usable in downstream tasks such as forecasting and control. To address this issue, we develop Constrained DEIM (C-DEIM) whose reconstructions are guaranteed to respect the physical bounds of the quantity of interest. C-DEIM enforces the bounds as soft constraints, in the form of a carefully designed penalty term, added to the underlying least squares problem. We prove that the C-DEIM reconstructions satisfy the physical constraints asymptotically, i.e., as the penalty parameter increases towards infinity. We also derive a quantitative upper bound for the observation residual of C-DEIM. Based on these theoretical results, we devise an efficient algorithm for practical implementation of C-DEIM. The efficacy of the method and the accompanying algorithm are demonstrated on several examples, including a heat transfer problem from fluid dynamics and a cellular automaton model of wildfire spread.
\end{abstract}


\end{frontmatter}

\section{Introduction}\label{sec:intro}
Discrete Empirical Interpolation Method (DEIM) is a framework for estimating a function from its partial pointwise observations. Such estimations are routinely needed for downstream tasks such as forecasting \cite{ghadami2022}, control \cite{da2021}, classification \cite{ahmed2021}, and sensitivity analysis \cite{hoffman2020}. 
DEIM is a linear method in that it estimates the unknown function as a linear combination of empirically prescribed basis functions~\cite{barrault2004,chaturantabut2010,drmac2016}. As such, it is simple, interpretable, and computationally inexpensive, in comparison to other reconstruction methods, such as autoencoders~\cite{Dawson2024}.

However, DEIM may produce physically inconsistent reconstructions when applied to a function with a restricted range. In many physical systems, quantities of interest have naturally imposed upper and/or lower bounds. For instance, mass density is intrinsically nonnegative. Likewise, quantities describing the probability of an event must lie between 0 and 1. In such cases, DEIM may generate reconstructions with values outside the physical bounds, leading to substantial reconstruction errors. Moreover, these unphysical values may impede downstream tasks such as forecasting, control, and uncertainty quantification.

To address this issue, here we propose Constrained DEIM (C-DEIM) which incorporates the constraints directly into the reconstruction process to avoid violations of known upper and lower bounds. We prove the convergence of the resulting C-DEIM algorithm and derive a rigorous upper bound for its observation residual. The contributions of this paper can be summarized as follows.
\begin{enumerate}
	\item We develop the C-DEIM method by incorporating the lower and upper bounds as a judiciously designed soft constraint in the underlying least squares problem.
	\item Convergence: We prove that, as the penalty parameter in the soft constraint increases, the C-DEIM solution converges to a reconstruction respecting the constraints.
	\item Error upper bound: We derive a rigorous upper bound for the observation residual corresponding to the C-DEIM reconstruction.
	\item C-DEIM algorithm: We propose an efficient implementation of the C-DEIM method. It combines Newton iterations with a bisection method, for determining the optimal value of the penalty parameter.
	\item Demonstration on numerical examples: We use three numerical examples, with increasing level of complexity, to demonstrate the efficacy and practicability of C-DEIM. In these examples, both optimal and randomized sensor placement methods are explored.
\end{enumerate}
 
\subsection{Related work}
Barrault et al.~\cite{barrault2004} introduced Empirical Interpolation Method (EIM) for reduced-order modeling of partial differential equations. Later, Sorensen and Chaturantabut~\cite{chaturantabut2010} developed the finite-dimensional version of EIM, known as Discrete Empirical Interpolation Method (DEIM). DEIM is in principle similar to Proper Orthogonal Decomposition (POD) in the sense that it projects the full-order dynamics onto a lower dimensional linear subspace. But unlike POD which uses an orthogonal projection, DEIM involves a carefully chosen oblique projection. This oblique projection significantly reduces the computational cost of evaluating the nonlinear terms in the differential equation.

Manohar et al.~\cite{manohar2018} demonstrated that DEIM can also be used to approximate a function from its discrete values. This is particularly useful for reconstructing a spatially distributed physical quantity from partial observational data. In practice, the observations are obtained from a set of sensors which make measurements of the function values. 

It is well-known that the quality of the reconstructions depends significantly on the sensor locations~\cite{maday2008}.  Drma\v{c} and Gugercin~\cite{drmac2016} introduced the Q-DEIM method which uses the QR factorization with column pivoting in order to estimate the optimal placement of the sensors.
Implementation of Q-DEIM is straightforward since column-pivoted QR factorization is widely available in standard computing libraries. Furthermore, Q-DEIM offers a tighter theoretical error bound compared to the original selection strategy introduced in~\cite{chaturantabut2010}. 

The original formulation of DEIM assumed that the number of basis functions (or modes), denoted by $m$, equals the number of sparse measurements (sensors), denoted by $r$~\cite{chaturantabut2010,drmac2016}. Subsequent work has explored both the underdetermined case ($r<m$) and the overdetermined case ($r>m$). In the overdetermined regime ($r>m$), where more sensors than modes are used, DEIM reconstructions show improved stability to observational noise compared to when the number of sensors is equal to the number of modes ($r=m$)~\cite{argaud2017, peherstorfer2020}. Likewise, Clark et al.~\cite{clark2020} demonstrated that increasing the number of sensors ($r>m$), for a fixed number of modes, leads to reduced reconstruction errors. The overdetermined case is equivalent to the gappy POD algorithm developed by Everson and Sirovich~\cite{everson1995}, a widely used technique for field reconstruction from partial observations~\cite{bui2004,willcox2006}. 

Relying on more sensors than modes becomes impractical if the sensors are expensive or if certain regions of the spatial domain are inaccessible~\cite{clark2020}. This motivated the investigation of the underdetermined regime ($r<m)$. In this regime, the (Q-)DEIM reconstruction error can be large~\cite{clark2020}. To improve reconstruction accuracy in the underdetermined case, Farazmand~\cite{farazmand2024} introduced Sparse DEIM (S-DEIM), which augments DEIM's minimum-norm solution with an optimally selected kernel vector. The resulting S-DEIM estimates have proved significantly more accurate than the (Q)-DEIM reconstructions~\cite{Farazmand2025}.

Other notable contributions include Extended DEIM (E-DEIM) which proposes a novel sensor placement method in the overdetermined case ($r>m$)~\cite{hendryx2021}. 
Localized DEIM (L-DEIM) approximates the function in a locally defined empirical basis instead of a global subspace~\cite{Willcox2014}.  
In the context of dynamical systems, Generalized DEIM (G-DEIM) identifies an optimal approximating subspace by leveraging the Jacobian information in addition to the state variables~\cite{Esmaeili2020}.
DEIM has also been extended for use with tensors~\cite{farazmand2023,kirsten2022}, weighted inner product spaces~\cite{drmac2018}, and for classification purposes~\cite{hendryx2021}.

Despite the numerous contributions reviewed above, existing DEIM-based methods do not explicitly account for range constraints on the function being interpolated. As mention earlier, such constraints arise routinely in physics and engineering. Interpolation and regression methods should respect these constraints in order for their outputs to be usable for downstream tasks such as forecasting and control.
Here we introduce C-DEIM to fulfill this requirement.

\subsection{Outline}
This paper is structured in the following manner. Section~\ref{sec:prelim} introduces the problem set-up and reviews DEIM. In Section~\ref{sec:cdeim}, we introduce C-DEIM and discuss its theoretical properties, including convergence analysis and error estimates. We develop the computational implementation of C-DEIM in Section~\ref{sec:prac_imp} and demonstrate its application on three numerical examples in Section~\ref{sec:numeric}. Finally, Section~\ref{sec:concl} offers our concluding remarks and outlines potential future work. 

\section{Preliminaries}\label{sec:prelim}
Consider a function $u:\Omega\to \mathbb R$ defined on a bounded spatial domain $\Omega\subset \mathbb R^d$.
We discretize this function over a fine spatial grid $\mathcal G\subset \Omega$ containing $N$ grid points. 
We assume that the function is known on a relatively sparse subset $\mathcal G_{s}\subset\mathcal G$. In practice, these known values are obtained from a set of $r$ sensors placed at $\mathcal G_s$. Our goal is to infer the function values over the entire computational domain $\mathcal G$.

It is convenient to express the discretized function as a vector $\vc u\in\mathbb R^N$.  The corresponding measurements at the sensor locations $\mathcal G_s$ are then given by a vector $\vc y \in \mathbb R^r$ which is a subset of the entries of $\vc u$. In practice, the observations $\vc y+\pmb\varepsilon$ contain measurement noise $\pmb\varepsilon$, which can be treated as a random variable in $\mathbb R^r$. In this section, for simplicity of the exposition, we neglect this noise, setting $\pmb\varepsilon=\mathbf 0$. In Section~\ref{sec:numeric}, where we present our numerical results, we further comment on the effect of the observational noise.
Defining the \emph{selection matrix} $C\in\mathbb R^{r\times N}$, which is an appropriate subset of the rows of the $N\times N$ identity matrix, we have $\vc y = C\vc u$. We further discuss the choice of the selection matrix $C$ in Section~\ref{sec:sensor_placement}.

We estimate the function $\vc u$ within an orthonormal set $\left\{\bphi_1,\cdots,\bphi_m\right\}\subset \mathbb R^N$ so that 
\begin{equation}
	\vc u \simeq \sum_{i=1}^m \alpha_i \bphi_i=\Phi\balpha.
\end{equation}
where $\Phi = [\bphi_1|\bphi_2|\cdots|\bphi_m]\in\mathbb R^{N\times m}$ is the basis matrix and $\balpha = [\alpha_1,\alpha_2,\cdots,\alpha_m]^\top\in\mathbb R^m$ are the coordinates in that basis. In principle, any appropriate set of orthonormal bases can be used; however, in practice, usually POD is used to extract a set of bases, i.e., POD modes~\cite{lumey2012,sirovich1987}.

Defining the total approximation error, 
\begin{equation}
\hat{\err}_m(\balpha) = |\Phi\balpha - \vc u|,
\label{eq:err_total}
\end{equation}
the optimal coefficients $\balpha$ that minimize the error $\hat{\err}_m$ are given by $\balpha  = \Phi^\top \vc u$, resulting in the corresponding \emph{best fit} reconstruction,
\begin{equation}
\hat{\vc u}=\Phi \Phi^\top \vc u.
\label{eq:best_fit}
\end{equation}
This is the best estimate of the vector $\vc u$ within the range of $\Phi$, spanned by the vectors $\{\bphi_1,\cdots,\bphi_m\}$.
The best fit reconstruction is not computable from sparse sensor measurements, since $\hat{\vc u}$ relies on the knowledge of the full state vector $\vc u$.

DEIM minimizes the discrepancy between the reconstructed function at the sensor locations and the true observations $\vc y=C\vc u$. This discrepancy is defined by the observation residual, 
\begin{equation}
\hat{\err}_o(\balpha)  =|C\Phi\balpha -\vc y|,
\label{eq:err_obs}
\end{equation}
where $|\cdot|$ denotes the standard Euclidean norm.
Note that $C\Phi\balpha$ corresponds to the entries of the reconstruction $\Phi\balpha$ at the sensor locations.

To obtain the optimal coefficients $\balpha$, we seek to minimize the observation residual $\hat{\err}_o(\balpha)$. Without loss of generality, it is more convenient to solve the convex minimization problem, 
\begin{equation}
\balpha_D=\argmin_{\balpha \in \mathbb{R}^m}\err_o(\balpha),
\label{eq:deim_opt}
\end{equation}
where $\err_o(\balpha):=\frac{1}{2}\hat{\err_o}(\balpha)^2$. Note that~\eqref{eq:deim_opt} is a linear least-squares problem whose minimum-norm solution is given by $\balpha_{D} = (C\Phi)^+\vc y$, where the superscript $+$ denotes the Moore-Penrose pseudo-inverse. The corresponding DEIM reconstruction is given by
\begin{equation}
\tilde{\vc u}_D = \Phi(C\Phi)^+\vc y.
\end{equation}
Note that the DEIM reconstruction $\tilde{\vc u}_D$ is computable because it relies on the observations $\vc y$ and not the entire vector $\vc u$.  We refer to the matrix $\Theta:=C\Phi$ as the \emph{sampled basis matrix} since its rows are a subset of the rows of $\Phi$, corresponding to the sensor locations.
Table~\ref{tab:symbols} summarizes the key mathematical quantities and notation used in this paper.
\begin{table}[ht]
	\centering
	\begin{tabular}{|c|l|}
		\hline
		\textbf{Symbol} & \textbf{Description} \\
		\hline \hline
		$N$ & High resolution grid size \\ \hline
		$m$ & Number of basis vectors \\ \hline
		$r$ & Number of sensors \\ \hline
		$\lambda$ & Penalty parameter \\ \hline
		$\mathcal P$ & Total penalty \\ \hline
		$\vc u \in \mathbb{R}^N$ & Discretized function \\ \hline
		$\tilde{\vc u} \in \mathbb{R}^N$ & Discretized reconstruction\\ \hline
		$\Phi \in \mathbb{R}^{N\times m}$ & Basis matrix \\ \hline
		$\balpha \in \mathbb{R}^m$ & Expansion coefficients, $\tilde{\vc u}=\Phi \balpha$ \\ \hline
		$C \in \mathbb{R}^{r \times N}$ & Selection matrix \\ \hline
		$\vc y \in \mathbb{R}^r$ & Observations, $\vc y=C\vc u$\\ \hline
		$\Theta \in \mathbb{R}^ {r \times m}$ & Sampled basis matrix, $\Theta=C\Phi$ \\
		\hline
	\end{tabular}
	\caption{The variables used in this paper.}
	\label{tab:symbols}
\end{table}

If $C$ is equal to the identity matrix, then the DEIM reconstruction is equivalent to the best fit reconstruction~\eqref{eq:best_fit}. However, here we only have access to sparse observations $\vc y=C\vc u$, where $C$ is a subset of the rows of the identity matrix; see Section~\ref{sec:sensor_placement}.

\subsection{Sensor Placement}\label{sec:sensor_placement}
The quality of the DEIM reconstruction is known to depend significantly on sensor locations $\mathcal G_s$ which are encoded in the selection matrix $C$. The optimal sensor placement minimizes the discrepancy between the DEIM reconstruction $\tilde{\vc u}_D$ and the true state $\vc u$. However, even if $\vc u$ was known, searching across all possible sensor configurations would be computationally intractable. Therefore, a greedy approach based on column-pivoted QR (CPQR) factorization is widely used~\cite{drmac2016}. Here, we briefly review this sensor placement method. Furthermore, towards the end of this section, we extend CPQR to the case where some subset of the spatial domain is inaccessible to sensors.

The rank-revealing QR decomposition of $\Phi^\top$ with column pivoting satisfies
\begin{equation}\label{eq:CPQR}
\Phi^\top \Pi = QR,
\end{equation}
where $\Pi \in \mathbb R^{N\times N}$ is a permutation matrix and $Q\in\mathbb R^{m\times m}$ is an orthogonal matrix. Then the Q-DEIM selection matrix is given by $C = \Pi_r^\top$ where $\Pi_r$ denotes the first $r$ columns of the permutation matrix $\Pi$. Roughly speaking, this Q-DEIM sensor placement minimizes a theoretical upper bound on the DEIM reconstruction error. We refer to Ref.~\cite{drmac2016} for a detailed discussion.

Q-DEIM selects $r$ sensor locations from the entire spatial grid $\mathcal{G}$. However, in many practical applications, feasible sensor locations are limited since some regions of the spatial domain may be inaccessible or too costly to monitor. Here, we propose a straightforward modification of CPQR which makes it applicable when feasible sensor locations are restricted.

Let $\mathcal G_{a}\subset \mathcal{G}$ represent the subset of the spatial grid which is accessible for sensor placement. In order words, $\mathcal G_a$ denotes the candidate sensor locations.
Note that the Q-DEIM sensor locations $\mathcal{G}_s$ may not necessarily lie within the accessible locations $\mathcal G_{a}$, i.e., $\mathcal{G}_s\not\subset\mathcal G_{a}$. To remedy this issue, we define the accessible selection matrix $C_a$ which is the $N\times N$ identity matrix with rows corresponding to inaccessible locations $\mathcal G\backslash \mathcal G_a$ zeroed out.
We then apply CPQR to $(C_a\Phi)^\top$, 
\begin{equation}\label{eq:rCPQR}
	(C_a\Phi)^\top \Pi = QR,
\end{equation}
which we refer to as the \emph{restricted CPQR} decomposition. Similar to Q-DEIM, we set the selection matrix $C=\Pi_r^\top$, where $\Pi_r$ again denotes the first $r$ columns of the permutation matrix $\Pi$. Since rows of $C_a\Phi$ corresponding to inaccessible grid points are zeroed out, the inaccessible region is demoted in the selection process.
This ensures that the selected sensor locations belong to the accessible set $\mathcal{G}_a$ while still aiming to minimize the reconstruction error. 

When the selection matrix $C$ is obtained from the CPQR decomposition~\eqref{eq:CPQR}, the sampled basis matrix $\Theta = C\Phi$ has full rank~\cite{Eisenstat1996}.
However, when the selection matrix is obtained from the restricted CPQR~\eqref{eq:rCPQR}, the resulting $\Theta$ is not guaranteed to have full rank. Nonetheless, in practice, this full-rank property tends to hold even if restricted CPQR is used for sensor placement. As a result, we make the following assumption when necessary.
\begin{ass} \label{ass:full_rank}
For a given basis matrix $\Phi$ and a selection matrix $C$, we assume the sampled basis matrix $\Theta = C\Phi\in\mathbb{R}^{r \times m}$ has full rank.
\end{ass}

\section{Constrained DEIM}{\label{sec:cdeim}}
While (Q-)DEIM reconstructs a field from sparse measurements, it does not account for the quantity's known physical constraints. In many applications, the quantity must lie within a restricted range. For example, the concentration of a substance is always a non-negative quantity. Here, we consider the general case, where the 
quantity of interest lies within an interval $[u_{\min},u_{\max}]$. Following the notation introduced in Section~\ref{sec:prelim}, the discretized quantity $\vc u = [u_1,u_2,\cdots, u_N]^\top$ satisfies
\begin{equation}\label{eq:bounds}
u_{\min} \leq u_i \leq u_{\max},\quad \text{for all}\quad i \in \{1, 2, \ldots, N\}.
\end{equation}
However, (Q-)DEIM offers no guarantee that the reconstructed function $\tilde{\vc u}_D$ will satisfy these bounds.
To address this issue, we propose a constrained extension of DEIM, called Constrained DEIM, which ensures that the reconstruction respects the prescribed upper and lower bounds~\eqref{eq:bounds}.

More specifically, in C-DEIM, we seek the optimal expansion coefficients $\balpha$ by solving the constrained optimization problem,
\begin{subequations}\label{eq:c_opt_1}
\begin{equation}
\min_{\balpha \in \mathbb R^m} \err_o(\balpha):=\frac12|C\Phi\balpha -\vc y|^2,
\end{equation}
\begin{equation}\label{eq:c_opt_const}
	\text{subject to} \quad u_{\min} \leq [\Phi\balpha]_i \leq u_{\max}, \quad i=1, \hdots, N,
\end{equation}
\end{subequations}
where $[\Phi \balpha]_i$ denotes the $i$-th coordinate of the reconstruction $\tilde{\vc u} = \Phi\balpha$.
Note that this optimization problem is similar to its DEIM counterpart~\eqref{eq:deim_opt}, but here we add the bound constraints on the reconstruction.

We point out that a similar constrained optimization problem was introduced in~\cite{argaud2017}. However, the type and purpose of their constraints differ from those used here in~\eqref{eq:c_opt_const}. In \cite{argaud2017}, the constraints are imposed directly on the coefficients $\balpha$ to improve the stability of the reconstruction under observational noise. In contrast, C-DEIM imposes inequality constraints on each component of the reconstruction $\tilde{\vc u}$ to enforce the physical upper and lower limits of the true function $\vc u$.

Let $\mathcal{S}\in \mathbb{R}^m$ denote the set of all coefficient vectors $\balpha\in \mathbb{R}^m$ such that the resulting reconstruction satisfies the constraints, 
\begin{equation}
\mathcal{S}:=\{\balpha \in \mathbb R^m\quad |u_{\min} \leq [\Phi\balpha]_i \leq u_{\max}, \quad i=1, \hdots, N\}.
\label{eq:s}
\end{equation} 
Then the constrained optimization problem~\eqref{eq:c_opt_1} can be written equivalently as
\begin{equation}
\balpha_{\min}=\argmin_{\balpha\in \mathcal{S}}\err_o(\balpha).
\label{eq:c_opt}
\end{equation}
In practice, if the basis matrix $\Phi$ is chosen appropriately, so that its range closely approximates the states $\vc u$, the constraint set $\mathcal S$ will be nonempty. For instance, if the basis matrix is obtained from POD modes, increasing the number of modes will eventually result in a nonempty constraint set.

Equation~\eqref{eq:c_opt_1} is a convex optimization problem subject to $2N$ linear inequality constraints. Although theoretically well-posed, it does not have a closed-form solution. 
There exist robust numerical methods, such as interior-point methods~\cite{Wright2000}, for solving such optimization problems.
However, these methods can be computationally expensive, particularly when the number of grid points $N$ is large.
Furthermore, in our experience, these methods may fail to converge to a minimizer when the parameter space $\mathbb R^m$ is high-dimensional.

Here, we pursue an alternative optimization approach by treating the hard constraints~\eqref{eq:c_opt_const} as soft constraints. To this end, we minimize the penalized cost function,
 \begin{equation}
f_\lambda(\balpha):=\err_o(\balpha)+\lambda \mathcal{P}(\balpha),
\label{eq:obj_f}
\end{equation}
where $\lambda\geq 0$ is a penalty parameter and $\mathcal{P}:\mathbb R^m\to\mathbb R^+$ penalizes deviations from the constraints~\eqref{eq:c_opt_const}. Note that the penalty term $\mathcal{P}$ can be viewed as a regularizer which penalizes overfitting to the observations $\vc y$. 
We emphasize that the penalty parameter $\lambda$ should not be confused with a hyperparameter; its value changes systematically during the optimization, as detailed in Section~\ref{subsec:update}.

In the remainder of this section, we discuss the theoretical properties of our soft-constraint approach. Specifically, in Section~\ref{subsec:penalty}, we 
discuss the class of admissible penalty functions $\mathcal P$ and provide an explicit example from this class. In Section~\ref{subsec:concept_alg}, we discuss the convergence properties of the C-DEIM algorithm; in particular, we show that as the penalty parameter $\lambda$ increases, the minimizer of $f_\lambda$ converges to the boundary of the constraint set $\mathcal S$.
In Section~\ref{subsec:obs_err}, we derive an upper bound on the observation residual $\err_o$ of the C-DEIM reconstruction.

\subsection{Penalty Function}\label{subsec:penalty}
The penalty term $\mathcal P$ is defined such that the C-DEIM reconstruction $\tilde{\vc u} = \Phi \balpha$ satisfies the physical constraints $u_{\min}\leq\tilde{u}_i\leq u_{\max}$. Here, we define $\mathcal P$ in terms of a scalar nonnegative function $p:\mathbb R\to\mathbb R$ which applies componentwise to the reconstruction, $\mathcal P(\balpha) = \sum_i p(\tilde u_i)$.

To ensure well-posedness and convergence of the forthcoming C-DEIM algorithm, we require that the function $p$ satisfies the following properties.
\begin{ass}[Admissible Penalty Functions]\label{ass:penalty}
The penalty function $p: \mathbb R\to\mathbb R$ possesses the following properties:
\begin{enumerate}
\item $p(u)$ is nonnegative, convex, twice continuously differentiable, with a Lipschitz continuous second derivative $p''$. Furthermore, $p''(u)>0$ for all $u\notin[u_{\min},u_{\max}]$.
\item $p(u)=0$ for all $u\in[u_{\min}, u_{\max}]$, so reconstructions that satisfy the constraints do not incur a penalty.
\item $p(u)$ is monotonically decreasing for $u<u_{\min}$ and monotonically increasing for $u>u_{\max}$.
\end{enumerate}
\end{ass}

\begin{figure} 
	\centering
	\includegraphics[width=0.5\textwidth]{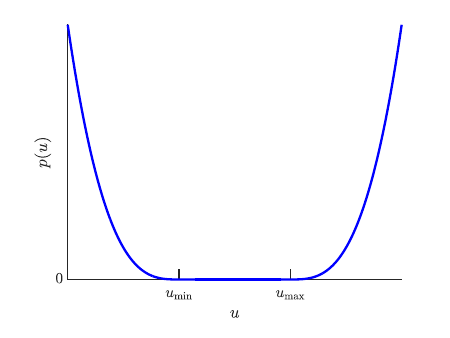}
	\caption{The piecewise cubic penalty function~\eqref{eq:p_fun}.}
	\label{fig:penalty}
	\end{figure}

An explicit example of a penalty function that satisfies the conditions of Assumption~\ref{ass:penalty} is given by the piecewise cubic function,
 \begin{equation}
	p(u)=\begin{cases}
		-\frac{1}{6}(u-u_{\min})^3 , \quad u<u_{\min},\\
		0 ,\quad u_{\min} \leq u \leq u_{\max},\\
		\frac{1}{6}(u-u_{\max})^3 , \quad u>u_{\max}.
	\end{cases}
	\label{eq:p_fun}
\end{equation}
The penalty function $p(u)$ is shown in figure \ref{fig:penalty}. Although we use this penalty function throughout our numerical examples (Section~\ref{sec:numeric}), 
any function that satisfies the conditions of Assumption~\ref{ass:penalty} can be used alternatively. For instance, C-DEIM with a piecewise fourth-order polynomial left the numerical reconstructions virtually unchanged, showcasing the robustness of the method to the choice of the penalty function.

We apply this penalty individually to each component of the reconstructed function. Using the Euclidean inner product, $\tilde u_i = \langle \Phi \balpha, \mathbf{e}_i\rangle$, where $\mathbf{e}_i$ is the $i$-th standard basis vector in $\mathbb{R}^N$, the total penalty is then defined as 
 \begin{equation}
 \mathcal{P}(\balpha)=\sum_{i=1}^N p(\langle \Phi \balpha, \mathbf{e}_i\rangle).
 \label{eq:penalty}
 \end{equation}
Note that $\mathcal P:\mathbb R^m\to\mathbb R^+$ is a convex, nonnegative function with Lipschitz-continuous second derivatives.
Consequently, the penalized cost function \eqref{eq:obj_f} is similarly convex and $C^2$ smooth. Furthermore, the set of coefficients for which the penalty $\mathcal P$ vanishes coincides with the constraint set~\eqref{eq:s}, i.e., $\mathcal S = \{\balpha\in\mathbb R^m: \mathcal P(\balpha) = 0 \}$. The penalty is identically zero inside the constraint set, which implies $\nabla\mathcal P(\balpha)=\vc 0$ for all $\balpha\in\mathcal S$. By continuity, the gradient of the penalty also vanishes on the boundary of the constraint set.

Although the total penalty~\eqref{eq:penalty} is the uniform sum of componentwise penalties, it can be easily generalized to a weighted sum where constraint violations are more severely penalized in certain spatial regions than others.
  
\subsection{A Conceptual C-DEIM Algorithm and its Convergence} \label{subsec:concept_alg}
Well-known results on convex optimization state that there exists a penalty parameter $\lambda>0$ such that a minimizer of the penalized cost function~\eqref{eq:obj_f} coincides with a minimizer of the constrained optimization problem~\eqref{eq:c_opt_1}; see, e.g.,~\cite{Lasserre2009}. However, these results do not generally specify the value of $\lambda$ for which this agreement occurs. Owing to the specific forms of the observation residual $\err_o(\balpha)$ and the penalty term $\mathcal P(\balpha)$, we can prove more specific results for C-DEIM. In particular, here we show that, as $\lambda$ tends to infinity, the minimizer of $f_\lambda$ converges to the boundary of the constraint set $\mathcal S$.

To state this result more concretely, we define
\begin{equation}
\balpha(\lambda):=\arg\min_{\balpha\in \mathbb R^m} f_\lambda(\balpha),
\label{eq:alpha_lam}
\end{equation}
as a minimizer of $f_\lambda(\balpha)$ for a given $\lambda>0$.

\begin{figure}
	\centering
	\includegraphics[width=.65\textwidth]{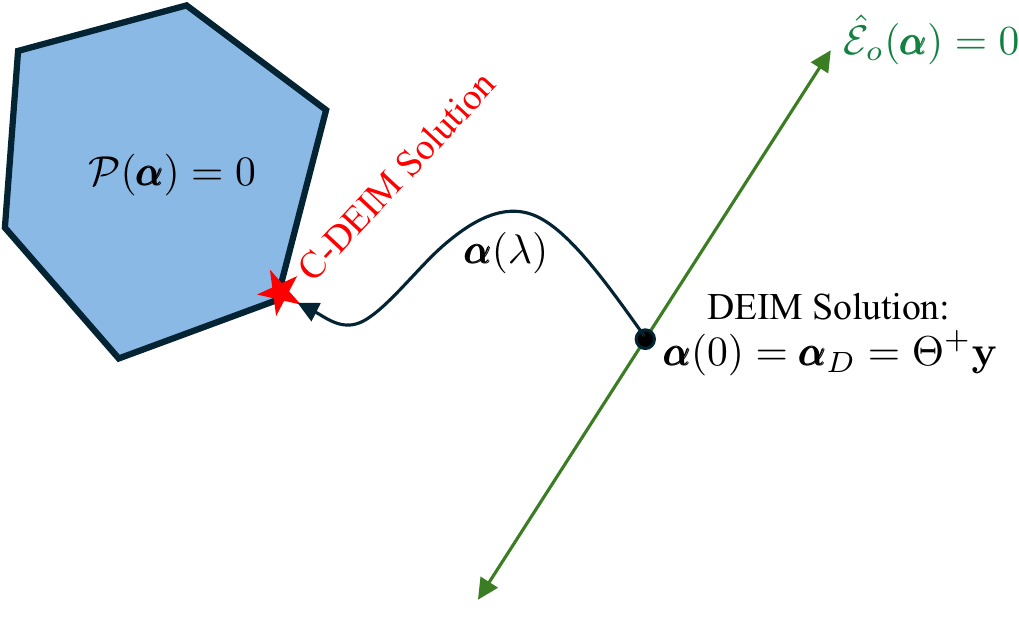}
	\caption{Schematic geometry of the C-DEIM method. As the penalty parameter $\lambda$ increases, the minimizer $\balpha(\lambda)$ converges to the constraint set $\mathcal S = \{\balpha\in\mathbb R^m: \mathcal P(\balpha) = 0 \}$. Note that the affine space, $\hat{\mathcal E}_o(\balpha)=0$, where the observation residual vanishes, may be an empty set in the overdetermined case ($r> m$).}
	\label{fig:constOpt}
\end{figure}
When $\lambda=0$, the penalty term vanishes, and the solution $\balpha(0)$ is equivalent to the DEIM solution, $\balpha(0)=\balpha_D=\Theta^+\vc y$. If this solution satisfies the constraints ($\balpha(0)\in\mathcal{S}$), the C-DEIM reconstruction uses the DEIM solution and no further action is needed. Otherwise, we seek a solution that satisfies the constraints by increasing $\lambda$ until the minimizer $\balpha(\lambda)$ falls within the constraint set $\mathcal{S}$; see figure~\ref{fig:constOpt}. The following theorem guarantees that this procedure eventually converges to a minimizer $\balpha(\lambda)$ that satisfies the constraints.
\begin{thm}[Convergence to the Constraint Set]
Consider a penalty function $\mathcal P(\balpha)$ satisfying the properties of Assumption~\ref{ass:penalty}. Then the minimizer~\eqref{eq:alpha_lam} satisfies,
\begin{equation}
\lim_{\lambda \to \infty} \id(\balpha(\lambda), \mathcal{S})=0,
\end{equation}
where $\id(\cdot,\cdot)$ denotes the Euclidean Hausdorff distance and $\mathcal S$ is the constraint set~\eqref{eq:s}.
\label{thm:convergence_lam}
\end{thm}
\begin{proof}
The proof is detailed in \ref{app:conv}.
\end{proof}

This theorem motivates the C-DEIM algorithm which is sketched conceptually in Algorithm~\ref{alg:cdeim_c}; a more detailed and practical algorithm is introduced in Section~\ref{sec:prac_imp}.
The algorithm iteratively increases the penalty parameter $\lambda$, starting from $\lambda=0$. At each iteration, it solves the minimization problem \eqref{eq:alpha_lam} using the current $\lambda$ value. This process continues until the resulting solution $\balpha(\lambda)$ generates a reconstruction $\tilde{\vc u}=\Phi\balpha(\lambda)$ that (approximately) satisfies the constraints $\mathcal{P}(\balpha(\lambda))=0$.
\begin{algorithm}
		\caption{Conceptual C-DEIM algorithm} \label{alg:cdeim_c}
		\textbf{Inputs:} \\
		Basis matrix $\Phi \in \mathbb{R}^{N \times m}$, sampled basis matrix $\Theta \in \mathbb{R}^{r \times m}$, observations $\vc y \in \mathbb{R}^{r}$
		\begin{algorithmic}[1]
			\State Initialize $\lambda=0$;
			\State $\balpha(0)=\Theta^+ \vc y$;
			\While{the penalty term $\mathcal{P}(\balpha(\lambda))>0$}
			\State Increase $\lambda$;
			\State Solve $\balpha(\lambda)=\arg\min_{\balpha} f_\lambda(\balpha)$;
			\EndWhile
			\State Compute reconstruction $\tilde{\vc{u}} = \Phi \balpha(\lambda)$.
				\end{algorithmic}
	\vspace{\baselineskip}
	\textbf{Output:} C-DEIM Reconstruction $\tilde{\vc{u}}\in \mathbb{R}^{N}$.
\end{algorithm}

Theorem~\ref{thm:convergence_lam} appears to imply that the penalty parameter $\lambda$ can be increased indefinitely. However, secondary considerations, such as the observation residual, limit the extent to which $\lambda$ can be increased in practice. We discuss this next.

\subsection{Observation Residual}\label{subsec:obs_err}
An attractive feature of DEIM is its interpolation property. If the number of sensors are less than or equal to the number of modes ($r\leq m$), then the DEIM reconstruction $\tilde{\vc u}_D$ agrees exactly with the sensor measurements, i.e., $C\tilde{\vc u}_D = \vc y$~\cite{chaturantabut2010,drmac2018}. Consequently, the corresponding observation residual $\hat{\err}_o(\balpha_D)$ vanishes. In contrast, the C-DEIM approach~\eqref{eq:obj_f} introduces a penalty term to enforce constraints, which may slightly compromise the interpolation property.

It is natural therefore to ask about the degree to which the C-DEIM reconstruction deviates from the observational data. To answer this question, here we derive an upper bound on the C-DEIM observation residual which quantifies the worst case scenario.
\begin{thm}[C-DEIM Observation Residual] \label{thm:obs_err}
Let $r\leq m$ and Assumption~\ref{ass:full_rank} hold. Then for any penalty parameter $\lambda>0$, the C-DEIM observation residual, corresponding to the minimizer~\eqref{eq:alpha_lam}, satisfies
\begin{equation}\label{eq:obs_ub}
\hat{\err}_o(\balpha(\lambda))\leq \frac{\lambda}{\sigma_{\min}(\Theta)} |\nabla \mathcal{P}(\balpha(\lambda))|,
\end{equation}
where $\sigma_{\min}(\Theta)$ is the smallest nonzero singular value of the sampled basis matrix $\Theta$.
\end{thm}
\begin{proof}
First note that since $\Theta$ has full row rank, $\Theta^+$ is its right inverse, 
\begin{equation}
		\Theta \Theta^+ = \left(\Theta^+\right)^\top \Theta^\top  = I_r,
		\label{eq:interp}
\end{equation}
where $I_r$ is the $r\times r$ identity matrix.
	
Next, consider the gradient of the penalized cost function~\eqref{eq:obj_f},
\begin{equation}
	\nabla f_\lambda(\balpha)=\Theta^\top\Theta\balpha-\Theta^\top\vc y +\lambda\nabla\mathcal{P}(\balpha).
\end{equation}
Since at the minimum $\nabla f_\lambda(\balpha(\lambda))=\mathbf{0}$, the minimizer must satisfy
	\begin{equation}
	\Theta \balpha(\lambda)=\vc y -\lambda\left(\Theta^+\right)^\top\nabla\mathcal{P}(\balpha(\lambda)).
	\end{equation}
	The corresponding C-DEIM reconstruction returns the following observation residual,
	\begin{equation}
	\hat{\err}_o(\balpha(\lambda))=|\Theta \balpha(\lambda) - \vc y| = \lambda|\left(\Theta^+\right)^\top\nabla\mathcal{P}(\balpha(\lambda))|. 
	\end{equation}
	Then, by definition of the 2-norm of a matrix and the fact that $\|A\|_2=\|A^\top\|_2$ for any matrix $A$, we obtain
	\begin{equation}
	\hat{\err}_o(\balpha(\lambda))\leq \lambda\|\Theta^{+}\|_2  |\nabla\mathcal{P}(\balpha(\lambda))|.
	\end{equation}
Finally, we obtain the desired result by recalling that $\|\Theta^+\|_2 = 1/\sigma_{\min}(\Theta)$.
\end{proof}

The upper bound~\eqref{eq:obs_ub} reveals the tension within C-DEIM between satisfying the constraints $\mathcal S$ and preserving the interpolation property. On one hand, Theorem~\ref{thm:convergence_lam} implies that larger penalty parameters $\lambda$ lead to better satisfaction of the constraints, so that $\lim_{\lambda\to\infty} |\nabla \mathcal P(\balpha(\lambda))|=0$.
On the other hand, the term $\lambda|\nabla \mathcal P(\balpha(\lambda))|$ in the upper bound~\eqref{eq:obs_ub} may converge to a nonzero value or even blow up as $\lambda$ increases, indicating that the observation residual can grow with $\lambda$. Consequently, the penalty parameter cannot be increased indefinitely without compromising the agreement between the reconstruction $\tilde{\vc u}$ and the observational data. This fact plays an important role in the numerical implementation of the C-DEIM algorithm; see Section~\ref{subsec:update}. 

We emphasize that Theorem~\ref{thm:obs_err} applies to the regime $r\leq m$, since the interpolation property~\eqref{eq:interp} does not generally hold in the overdetermined regime where $r>m$. This does not preclude the use of C-DEIM in the overdetermined regime ($r>m$); however, the derived error bound on the observation residual does not apply. Numerical results for the observation residual in the overdetermined regime are presented in figure~\ref{fig:cdeim_errors_modes}(c).

\section{Numerical Implementation of C-DEIM}\label{sec:prac_imp}
Theorem~\ref{thm:convergence_lam} established that the C-DEIM reconstruction $\balpha(\lambda)$ converges to the constraint set $\mathcal{S}$ as the penalty parameter $\lambda$ increases to infinity. To guarantee the C-DEIM algorithm halts in finite iterations, we introduce a small tolerance parameter $\delta>0$ and stop the algorithm when
\begin{equation}
\mathcal{P}(\balpha(\lambda))<\delta.
\label{eq:practical_stop}
\end{equation}
Since the total penalty $\mathcal{P}(\balpha(\lambda))$ is zero if and only if the reconstruction $\tilde{\vc u}=\Phi\balpha(\lambda)$ strictly satisfies the constraints, a small positive $\delta$ allows for slight, controlled violations. 

Using the penalty function \eqref{eq:p_fun} and a given $\delta >0$, the resulting C-DEIM reconstructions will satisfy
\begin{equation}
	u_{\min}-(6\delta)^{1/3} \leq \tilde{u}_i \leq u_{\max}+(6\delta)^{1/3},\quad i=1,2, \hdots, N.
	\label{eq:tol_bounds}
\end{equation}
This provides a guarantee that the reconstruction will be approximately within the prescribed bounds, with the maximum deviation controlled by the tolerance $\delta$.
If strict adherence to the prescribed bounds is required for downstream tasks, the reconstruction $\tilde{\vc u}$ can be post-processed by rounding any out-of-bounds values to $u_{\min}$ or $u_{\max}$.

C-DEIM requires the computation of the minimizer~\eqref{eq:alpha_lam} for a given penalty parameter $\lambda$. To estimate the minimizer $\balpha(\lambda)$, we use the standard Newton method, as detailed in Section~\ref{subsec:newton_iter}. Recall that unnecessarily large penalty parameters may result in large observation residuals; cf Section~\ref{subsec:obs_err}. Therefore, $\lambda$ needs to be chosen judiciously. Section~\ref{subsec:update} proposes a bisection method for determining the optimal penalty parameter $\lambda$.

\subsection{Newton Iterations}
\label{subsec:newton_iter}
For a fixed penalty parameter $\lambda$, C-DEIM requires minimizing the convex, penalized cost function~\eqref{eq:obj_f}. We use Newton's method to compute the minimizer $\balpha(\lambda)$, which satisfies $\nabla f_\lambda(\balpha(\lambda))=\vc 0$. The gradient and the Hessian of the penalized cost function \eqref{eq:obj_f} are given by
\begin{subequations}
\begin{equation}
\nabla f_\lambda(\balpha)= \Theta^{\top}\Theta\balpha-\Theta^\top \vc y+\lambda\Phi^\top p'(\Phi\balpha),
\end{equation}
\begin{equation}
\nabla^2f_\lambda(\balpha)=\Theta^{\top}\Theta+\lambda\Phi^\top D(\Phi\balpha)\Phi, \label{eq:hessian_f}
\end{equation}
\end{subequations}
where $p'$ is applied componentwise to the entries of the reconstruction $\tilde{\vc u}=\Phi\balpha$, and $D(\tilde{\vc u})$ is a diagonal matrix with diagonal entries $D_{ii}(\tilde{\vc u}) = p''(\tilde u_i)$. With a full-rank sampled basis matrix $\Theta$ (Assumption \ref{ass:full_rank}) and a penalty function $p$ that satisfies Assumption \ref{ass:penalty}, the Hessian \eqref{eq:hessian_f} is positive-definite and therefore nonsingular for all $\balpha(\lambda)\notin\mathcal{S}$ (See \ref{app:newton}).

Starting with an initial guess $\balpha_0$ for the minimizer, Newton's method iteratively updates $\balpha$ according to 
\begin{equation}
\balpha_{k+1}=\balpha_k-[\nabla^2 f_\lambda(\balpha_k)]^{-1}\nabla f_\lambda(\balpha_k),
\label{eq:newton_update} 
\end{equation}
where $\balpha_k(\lambda)$ is the estimate at the $k$-the iteration. For notational simplicity, we omit the dependence of the iterates $\balpha_k(\lambda)$ on the penalty parameter, and simply write $\balpha_k$. The Newton iterations~\eqref{eq:newton_update} continue until successive estimates are within a specified tolerance $\tau$, i.e., $|\balpha_k-\balpha_{k-1}|\leq \tau$. Newton's method is detailed in Algorithm~\ref{alg:newton_solv}.

\begin{algorithm}
		\caption{Newton solver} \label{alg:newton_solv}
		\textbf{Inputs:}
		Penalty parameter $\lambda$, tolerance $\tau$, initial guess $\balpha_{\text{init}}$
		\begin{algorithmic}[1]
				\State Initialize $\balpha_0=\balpha_{\text{init}}$, $k=1$;
				\While{$|\balpha_k-\balpha_{k-1}|>\tau$}
				\State $\balpha_{k+1}=\balpha_k-[\nabla^2 f_\lambda(\balpha_k)]^{-1}\nabla f_\lambda(\balpha_k$);\Comment{Compute Newton step.}
				\State $k=k+1$;
				\EndWhile
			\State $\balpha(\lambda)=\balpha_k$;\;
	\end{algorithmic}
	\vspace{\baselineskip}
	\textbf{Output:} Minimizer $\balpha(\lambda)$
	\end{algorithm}

To guarantee local convergence of Newton's method, the Hessian $\nabla^2f_\lambda(\balpha)$ must be nonsingular and Lipschitz continuous in a neighborhood of the minimizer~\cite{kelley2003}. Since the second derivative of the piecewise cubic penalty function \eqref{eq:p_fun} is piecewise linear, $p''$ is Lipschitz continuous. We show that the conditions of convergence are met for the function $f_\lambda(\balpha)$ under our assumptions in the following theorem.
\begin{thm}[Convergence of Newton Iterations]\label{thm:conv}
Let the penalty $\mathcal P$ satisfy Assumption~\ref{ass:penalty}.
Then for any $\lambda>0$ and $\balpha_0\notin\mathcal S$, the Newton iterations~\eqref{eq:newton_update} are locally convergent to the minimizer $\balpha(\lambda)$ of $f_\lambda$.
\end{thm}
\begin{proof}
See \ref{app:newton}.
\end{proof}
Although Theorem~\ref{thm:conv} guarantees the convergence of the Newton iterations for any $\lambda$, it does not provide a method for choosing this penalty parameter. In the next section, we discuss a bisection method for selecting the optimal value of the penalty parameter.

\subsection{Optimal Penalty Parameter}\label{subsec:update}
As discussed in Section~\ref{subsec:obs_err}, the penalty parameter $\lambda$ must balance the constraint enforcement against agreement with observations. A larger penalty parameter means better agreement with the constraints but may lead to a larger observation residual. Therefore, we seek the smallest possible $\lambda$ which satisfies the stopping criterion~\eqref{eq:practical_stop}. We denote this optimal value of the penalty parameter by $\lambda_{\text{opt}}$.

Our strategy to find this optimal $\lambda$ involves an initial broad search, followed by a bisection refinement. We start with $\lambda=0$. If the resulting solution $\balpha(0)$ satisfies the stopping criterion,  $\mathcal{P}(\balpha(0)) < \delta$, we are done. Otherwise, we find the minimizer $\balpha(\lambda)$ corresponding to $\lambda=\lambda_{\text{init}}$, where $\lambda_{\text{init}}$ is very small.
We successively increase $\lambda$ by a factor $\gamma>1$ until the stopping criterion is met for a value of $\lambda$, denoted by $\overline\lambda$. One can set $\gamma=2$ without affecting the forthcoming argument.
Let $\underline\lambda$ be the value of $\lambda$ in the previous step, $\underline\lambda=\overline\lambda/\gamma$,
satisfying $\mathcal{P}(\balpha(\underline\lambda)) \geq \delta$. This establishes an interval $[\underline\lambda,\overline\lambda]$ that encloses the desired optimal $\lambda_{\text{opt}}$.

Then, we use bisection within $[\underline\lambda,\overline\lambda]$ to refine the estimate of the optimal $\lambda$ that satisfies $\mathcal{P}(\balpha(\lambda))< \delta$. In each step, we test the midpoint of the interval, $\lambda_{\text{mid}}=(\underline\lambda+\overline\lambda)/2$. If $\mathcal{P}(\balpha(\lambda_{\text{mid}})) < \delta$, it means a smaller $\lambda$ will satisfy the criterion, so we redefine the upper bound $\overline\lambda=\lambda_{\text{mid}}$. Alternatively, if $\mathcal{P}(\balpha(\lambda_{\text{mid}})) \geq\delta$, the midpoint is not sufficiently large, so we increase the lower bound $\underline\lambda=\lambda_{\text{mid}}$. This continues until the width of the interval $[\underline\lambda,\overline\lambda]$ is less than a prescribed small bisection tolerance $\tau_{\lambda}$. Note that, by construction, $\lambda_{\text{opt}}\in[\underline\lambda,\overline\lambda]$ at every iteration of our bisection method.
The final value of $\overline\lambda$ is then taken as the optimal penalty parameter $\lambda_{\text{opt}}$.
\begin{algorithm}
	\caption{C-DEIM algorithm} \label{alg:practical}
	\textbf{Inputs:}
	\begin{itemize}
		\item Basis matrix, $\Phi \in \mathbb{R}^{N \times m}$, sampled basis matrix $\Theta \in \mathbb{R}^{r \times m}$, observations $\vc y \in \mathbb{R}^{r}$
		\item Initial penalty parameter $0<\lambda_{\text{init}}\ll 1$, growth factor $\gamma>1$
		\item Constraint tolerances $\delta$, bisection tolerance $\tau_{\lambda}$ 
	\end{itemize}
	\begin{algorithmic}[1]
		\State $\balpha=\Theta^+ \vc y$;
		\If{$\mathcal P(\balpha)<\delta$}
		\State $\lambda_{\text{opt}}=0$;
		\Else
		\State Initialize $\lambda=\lambda_{\text{init}}$;
		\While{$\mathcal{P}(\balpha)\geq\delta$}
		\State $\lambda=\gamma \lambda$; \Comment{Increase $\lambda$.}
		\State $\balpha=\texttt{NewtonSolve}(\balpha,\lambda)$; 
		\EndWhile
		\State $\underline\lambda=\lambda/\gamma$, $\overline\lambda=\lambda$; \Comment{Initialize bisection interval.}
		\While{$|\overline\lambda-\underline\lambda|>\tau_{\lambda}$}
		\State  $\lambda_{\text{mid}}=(\underline\lambda+\overline\lambda)/2$; \Comment{Compute midpoint.}
		\State $\balpha=\texttt{NewtonSolve}(\balpha,\lambda_{\text{mid}})$
		\If{$\mathcal{P}(\balpha)\geq\delta$} 
		\State Set $\underline\lambda=\lambda_{\text{mid}}$; \Comment{Raise lower bound.}
		\Else 
		\State Set $\overline\lambda=\lambda_{\text{mid}}$; \Comment{Lower upper bound.}
		\EndIf	
		\EndWhile
		\State $\lambda_{\text{opt}}=\overline\lambda$;
		\State $\balpha=\texttt{NewtonSolve}(\balpha,\lambda_{\text{opt}})$;
		\EndIf
		\State Compute the reconstruction $\tilde{\vc u}=\Phi \balpha$;		
	\end{algorithmic}
	\textbf{Output:} C-DEIM Reconstruction $ \tilde{\vc{u}}\in \mathbb{R}^{N}$
\end{algorithm} 

The bisection method described above ensures that the smallest penalty parameter satisfying the stopping criterion is selected. By avoiding excessively large $\lambda$ values, it keeps the observation residual $\hat{\err}_o(\balpha(\lambda))$ as small as possible while simultaneously providing a reconstruction that effectively respects the range constraints~\eqref{eq:c_opt_const}. 

Algorithm~\ref{alg:practical} summarizes the numerical implementation of C-DEIM, combining the adaptive penalty parameter update strategy (bisection method) with Newton's method~\eqref{eq:newton_update}. 
	
\section{Numerical Results}\label{sec:numeric}
We present three numerical examples to demonstrate the efficacy of C-DEIM. The first example is a benchmark problem detailed in Section~\ref{subsec:toy}, involving the interpolation of random harmonic functions. In the second example, we reconstruct the temperature field in the Rayleigh--B\'enard convection, which arises in fluid mechanics (Section~\ref{sec:RBC}).
In this example, the temperature of the fluid is bounded by its prescribed boundary conditions.

In Section~\ref{subsec:fire}, we consider cellular automata simulations of wildfire propagation. This model includes a state variable, representing the fraction of a cell that has been burned by the fire. As a result, the state variable is naturally restricted between 0 and 1. Using C-DEIM with sparse measurements of the fire, one hour after ignition, we generate a reconstruction. Furthermore, we demonstrate the utility of the physically consistent C-DEIM reconstruction by using it as an initial condition to forecast the future state of the fire. 

In each example, we compute the relative $L^2$ reconstruction error and relative observation residual, respectively, using
\begin{equation}
RE:=\frac{|\tilde{\vc u}-\vc u|}{|\vc u|}, \quad OR:=\frac{|C\tilde{\vc u}-\vc y|}{|\vc y|},
\end{equation}
where $|\cdot|$ denotes the Euclidean norm. Throughout this section, we set the number of basis functions equal to the number of available sensors ($r=m$), unless explicitly stated otherwise. 

We find that C-DEIM is remarkably robust to observational noise. In particular, the results are virtually unchanged when the observations are polluted with up to 10\% synthetic Gaussian noise. We attribute this robustness to the regularizing effect of the penalty term in the C-DEIM formulation~\eqref{eq:obj_f} which suppresses overfitting to observations. The results in this section are reported for noise-free observations.

\subsection{Benchmark: Random Harmonics}\label{subsec:toy}

We use a benchmark problem to demonstrate the results of C-DEIM. We consider a class of functions defined as a linear combination of cosine terms, 
	\begin{equation}
	g_j(x)=\sum_{k=1}^{20} a_{kj}\cos(kx+\phi_{kj}),
	\end{equation}
	where the random amplitudes $a_{kj}$ are drawn from a normal distribution $\mathcal{N}(0,1/k)$ and the random phases $\phi_{kj}$ are drawn from a uniform distribution $\mathcal{U}[0, 2\pi]$. To introduce known bounds, we normalize each function by its maximum absolute value, ensuring its range is within $[-1,1]$, 
	\begin{equation}
	 u_j(x)=\frac{ g_j(x)}{\max_{x} |g_j(x)|}.
	\end{equation}
Each function $u_j(x)$ is discretized over $N=1000$ equally spaced points on the interval $x\in[0,2\pi]$ to create a vector $\vc u_j$. We generate an ensemble of 1000 normalized functions. The first 800 realizations are used as training data; namely, to compute the basis matrix $\Phi$ using POD and to derive the selection matrix $C$ using CPQR.  Although here, the basis matrix $\Phi$ comprises POD modes, it could consist of any orthonormal set of basis functions, such as Fourier modes. We use the remaining 200 functions as test cases to compare the performance of C-DEIM and DEIM. 

When sensors are placed without restriction, Q-DEIM often produces reconstructions satisfying the constraints. However, this performance degrades when sensor placement is restricted, especially when sensors cannot be located near domain boundaries. In such cases, Q-DEIM reconstructions frequently violate the constraints, producing values outside the prescribed bounds. To model this scenario, we restrict the accessible sensor locations to the subset, $\mathcal G_a=[\eta,2\pi-\eta]$ with $\eta=0.1\pi$, and compare the reconstructions from both DEIM and C-DEIM. 

\begin{figure} 
	\centering
	\includegraphics[width=\textwidth]{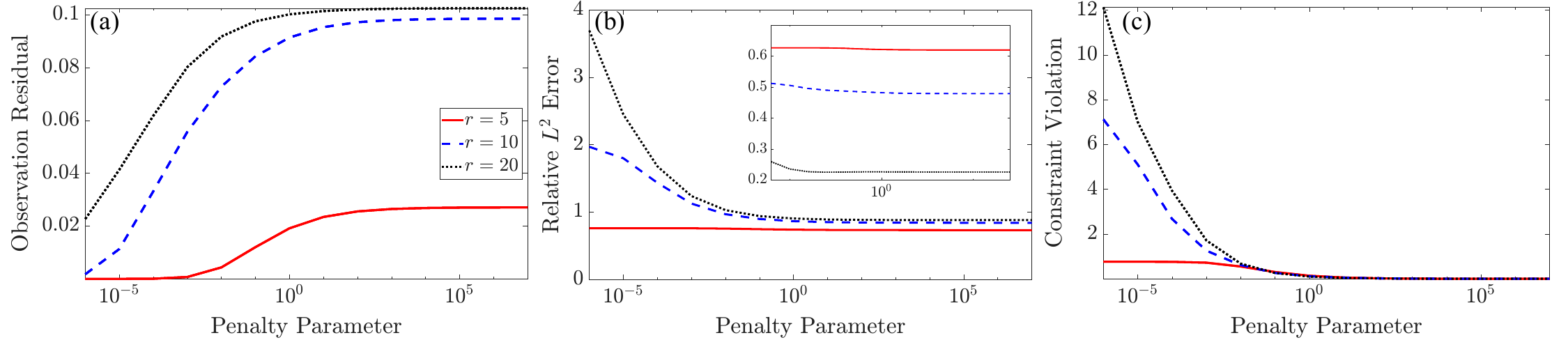}
	\caption{(a) Mean relative observation residual of C-DEIM reconstructions as a function of the penalty parameter $\lambda$, using $r$ sensors.  (b) Mean relative $L^2$ reconstruction error for C-DEIM as $\lambda$ increases. (inset) Mean relative $L^2$ reconstruction error over the accessible sensor region. (c) Maximum constraint violation for C-DEIM reconstructions as the penalty parameter $\lambda$ increases.}
	\label{fig:increase_lam}
	\end{figure}

In figure~\ref{fig:increase_lam}, we first discuss the behavior of the C-DEIM reconstruction with respect to the penalty parameter $\lambda$. Figure~\ref{fig:increase_lam}(a) shows that the mean relative observation residual increases with the penalty parameter $\lambda$, while figure~\ref{fig:increase_lam}(b) tracks the relative reconstruction error. This is in agreement with the theoretical results of Section~\ref{subsec:obs_err}. As shown in figure~\ref{fig:increase_lam}(c),  increasing $\lambda$ decreases the reconstruction error by reducing constraint violations, where the maximum constraint violation is defined as $\max \{ \max_i (\tilde u_i-u_{\max}), \max_i (u_{\min} - \tilde u_i)\}$, with $\tilde u_i$ denoting the components of the C-DEIM reconstruction $\tilde{\vc u}$. 
At the same time, increasing $\lambda$ also results in a higher observation residual. This observation residual eventually plateaus. This suggests an optimal $\lambda$ that balances these competing objectives. The goal of the bisection method introduced in Section~\ref{subsec:update} is to estimate this optimal $\lambda$. 

It is noteworthy that, although the mean relative reconstruction error decreases with larger penalty parameters, it remains high due to large errors in the extrapolation region, where there are no sensor measurements. As shown in the inset of figure~\ref{fig:increase_lam}(b), within the accessible sensor region $\mathcal G_a$, the mean relative reconstruction error decreases as the number of sensors increases. For $r=20$ sensors, this error plateaus around 22\% as the penalty parameter increases.

\begin{figure} 
	\centering
	\includegraphics[width=\textwidth]{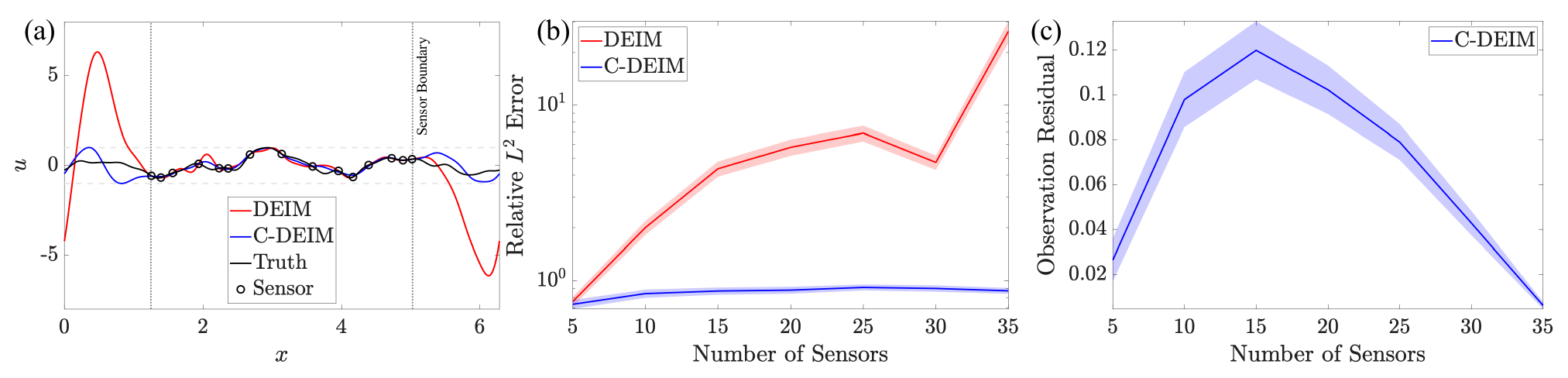}
	\caption{(a) Example comparison of truth (black), DEIM reconstruction (red) and C-DEIM reconstruction (blue) for one test snapshot with $r=15$ sensors. The DEIM reconstruction violates the bounds, while the C-DEIM reconstruction stays within $[-1,1]$. (b) Mean relative $L^2$ reconstruction error across 200 test snapshots for increasing number of sensors. C-DEIM achieves significantly lower reconstruction errors than DEIM, which frequently produces large reconstruction errors due to constraint violations. Shaded regions represent a 95\% confidence interval. (c) C-DEIM mean relative observation residual. Even with deviations from observed values, observation residual remains low. Shaded regions show the 95\% confidence interval.}
	\label{fig:cdeim_errors_all}
	\end{figure}

Figure~\ref{fig:cdeim_errors_all} shows overall performance metrics for the reconstruction of the test data. Figure~\ref{fig:cdeim_errors_all}(a) illustrates a single representative sample from the test data and its reconstruction, using $r=15$ sensors placed in the restricted subdomain $\mathcal G_a$. Inside this interpolation region, the DEIM and C-DEIM reconstructions are very close to each other and in good agreement with the true function. However, in the extrapolation region $\mathcal G\backslash \mathcal G_a$, the DEIM reconstruction deviates substantially from the bounds  $ -1\leq u_j(x)\leq 1$. In contrast, the C-DEIM reconstruction respects the constraints while achieving a low relative observation residual of 10.3\%. For this function, the optimal penalty parameter is $\lambda_{\text{opt}}=156.25$.

In this example, we stop the C-DEIM algorithm when the total penalty $\mathcal{P}(\balpha(\lambda))$ falls below the prescribed tolerance $\delta=10^{-7}$. This choice of $\delta$ guarantees that each component of the reconstruction can attain at most a small deviation from the prescribed bounds, i.e.  $\tilde{u}_i\in[-1.0084,1.0084]$; see equation~\eqref{eq:tol_bounds}. This demonstrates that the practical stopping criterion provides a strict control over the maximum constraint violation. We use an initial penalty parameter $\lambda_{\text{init}}=10^{-7}$, Newton tolerance $\tau=10^{-10}$ and growth factor $\gamma=10$.
We assess the performance of C-DEIM and DEIM in figures~\ref{fig:cdeim_errors_all}(b) and (c) across all 200 test snapshots. Figure~\ref{fig:cdeim_errors_all}(b) compares the mean relative reconstruction errors as the number of sensors increases. Both reconstructions have large relative errors because of large deviations in the extrapolation region $\mathcal G\backslash \mathcal G_a$. Nonetheless, the C-DEIM reconstructions return a significantly smaller relative error. The mean relative reconstruction errors remain below 100\% for C-DEIM. In contrast, the DEIM error grows as the number of sensors increases, reaching over 2000\% mean relative error for $r=35$ sensors. 

Figure~\ref{fig:cdeim_errors_all}(c) shows the mean relative observation residual $\hat{\mathcal E}_o$ for C-DEIM, which quantifies the deviation of the reconstruction from sensor measurements. Recall from Section~\ref{subsec:obs_err} that the DEIM reconstruction respects the interpolation property, and therefore it agrees exactly with the observations, $\hat{\err}_o(\balpha_D)=0$. On the other hand, the C-DEIM reconstructions can return a nonzero observation residual.
Figure~\ref{fig:cdeim_errors_all}(c) shows that the C-DEIM observation residual is nonetheless small. More specifically, the mean observation residual of C-DEIM, averaged across the test data, does not exceed 12\%.

Compared to DEIM, our C-DEIM method requires the additional steps of solving the Newton iterations \eqref{eq:newton_update} and determining the optimal penalty parameter $\lambda$ using the bisection method detailed in Section~\ref{subsec:update}. Even so, the C-DEIM algorithm converges relatively quickly. For example, to reconstruct one of the test snapshots, C-DEIM uses the DEIM solution $\balpha_D$ as its initial guess, and then, on average, requires an additional 0.27 seconds to arrive at the optimal coefficients $\balpha(\lambda_\text{opt})$. Therefore, despite requiring additional steps, the computational cost of C-DEIM remains comparable to DEIM. 

\begin{figure} 
	\centering
	\includegraphics[width=\textwidth]{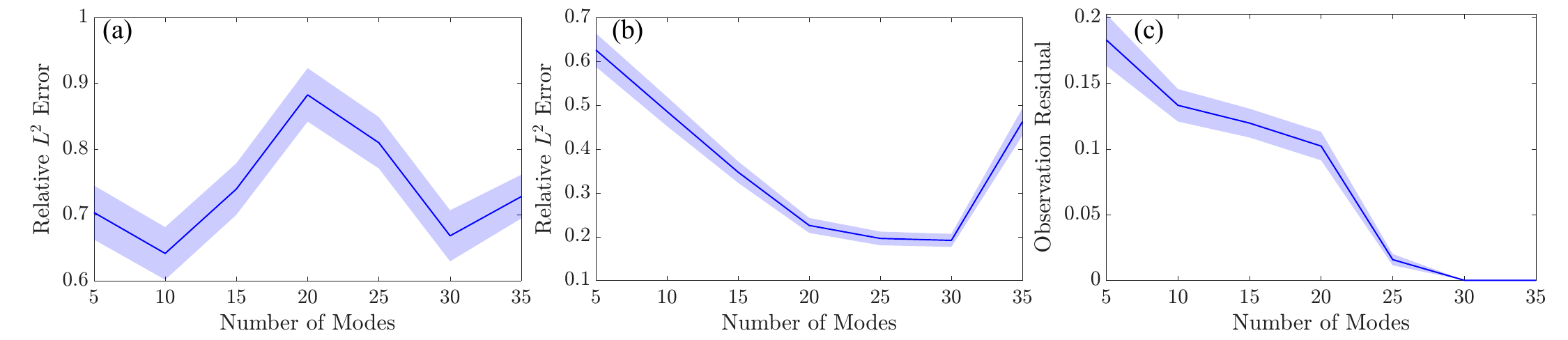}
	\caption{Mean relative reconstruction error with $r=20$ sensors and an increasing number of modes. (a) Mean relative reconstruction error across the entire domain. (b) Mean relative reconstruction error within the accessible sensor region. (c) Relative observation residual. Shaded regions represent a 95\% confidence interval.}
	\label{fig:cdeim_errors_modes}
\end{figure}

Following common practice~\cite{chaturantabut2010,drmac2016}, here we reported results for equal number of sensors and modes, $r=m$. For completeness, we also assess the performance of C-DEIM in the underdetermined ($r<m$) and overdetermined ($r>m$) regimes. To this end, we fix the number of sensors at $r=20$ and vary the number of modes $m$. As shown in figure~\ref{fig:cdeim_errors_modes}(a), the total mean relative error peaks when the number of sensors is equal to the number of modes at $m=r=20$. However, the relative reconstruction error within the accessible sensor region $\mathcal G_a$ decreases as the number of modes increases, until it increases again at $m=35$; see figure~\ref{fig:cdeim_errors_modes}(b). Meanwhile, the relative observation residual decreases monotonically as the number of modes increases (figure \ref{fig:cdeim_errors_modes}(c)), since the higher model complexity allows for better fits to the observations. We note that for $r<m$, the problem is underdetermined and the DEIM minimizer $\pmb\alpha_D$ is not unique; here, we report results corresponding to the unique minimum-norm solution.

\subsection{Rayleigh--B\'enard Convection} \label{sec:RBC}
As our next example, we consider the Rayleigh--B\'enard convection in a two-dimensional rectangular domain $(x_1,x_2)\in[0,L]\times [0,H]$. An incompressible fluid is confined in this domain and heated by the lower wall at $x_2=0$ and cooled by the upper wall at $x_2=H$. These horizontal walls are kept at a constant temperature. Denoting the fluid velocity field by $\vc v(\vc x,t)$ and the fluid temperature by $T(\vc x,t)$, the equations of motion, in dimensionless variables, are given by
\begin{subequations}\label{eq:rbc}
\begin{equation}\label{eq:rbc-v}
	\pard{\vc v}{t}+\vc v\cdot\nabla\vc v = -\nabla p + T\vc e_2+\sqrt{\frac{Pr}{Ra}}\Delta\vc v,\quad \nabla\cdot\vc v=0,
\end{equation}
\begin{equation}\label{eq:rbc-T}
	\pard{T}{t}+\vc v\cdot\nabla T = \frac{1}{\sqrt{PrRa}}\Delta T,
\end{equation}
\end{subequations}
where $p(\vc x,t)$ is the pressure field and $\vc e_2$ is the unit vector pointing in the vertical direction $x_2$~\cite{Kooloth2021}. Here, the Prandtl number $Pr$ and the Rayleigh number $Ra$ are dimensionless quantities. We numerically solve these equations with $Pr=10$ and $Ra=10^5$, on a rectangular domain with dimensionless lengths $L=4$ and $H=1$. The stress-free boundary conditions~\cite{Lohse2014} are given by
\begin{align}
v_1=0, \quad \pard{v_2}{x_1}=0,\quad \pard{T}{x_1}=0, & \quad \mbox{at}\quad	x_1 \in\{0, 4\}\ \mbox{(vertical walls)},\nonumber\\
\pard{v_1}{x_2}=0, \quad v_2=0,\quad T=1, &  \quad  \mbox{at}\quad x_2=0\ \mbox{(bottom wall)},\nonumber\\
\pard{v_1}{x_2}=0, \quad v_2=0,\quad T=0, &  \quad \mbox{at}\quad  x_2=1\ \mbox{(top wall)}.
\label{eq:RBC_bc}
\end{align}

\begin{figure}
	\centering
	\includegraphics[width=\textwidth]{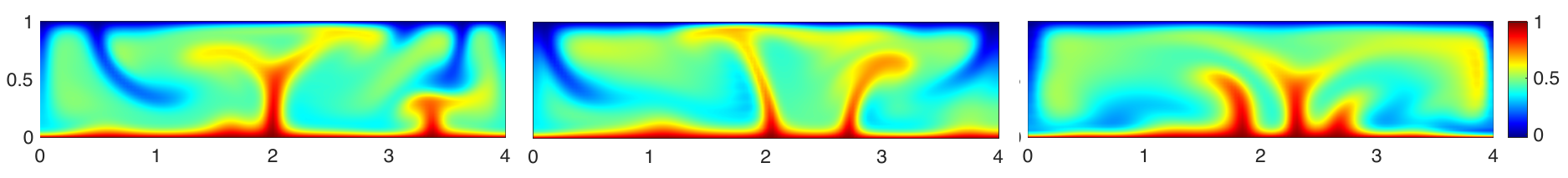}
	\caption{First three snapshots of the temperature from the training set of the Rayleigh--B\'enard convection. The snapshots are five time units apart.}
	\label{fig:RBC_snapshots}
\end{figure}

Note that the quiescent flow ($\vc v=0$) with the linear temperature profile, $T(\vc x,t) = 1-x_2$, is an exact solution of~\eqref{eq:rbc}. However, this solution is an unstable equilibrium. We initialize our simulation from a small perturbation to the quiescent flow. This generates a turbulent flow after a short transient excursion. After the initial transients have passed, we record the solution at intervals of $\Delta t=5$ time units, resulting in 1000 snapshots altogether. The first three snapshots are shown in figure~\ref{fig:RBC_snapshots}. We use the first 950 snapshots as training data, from which the POD modes and the sensor locations are computed. The remaining $50$ snapshots are used as test data.

\begin{figure}
	\centering
	\includegraphics[width=\textwidth]{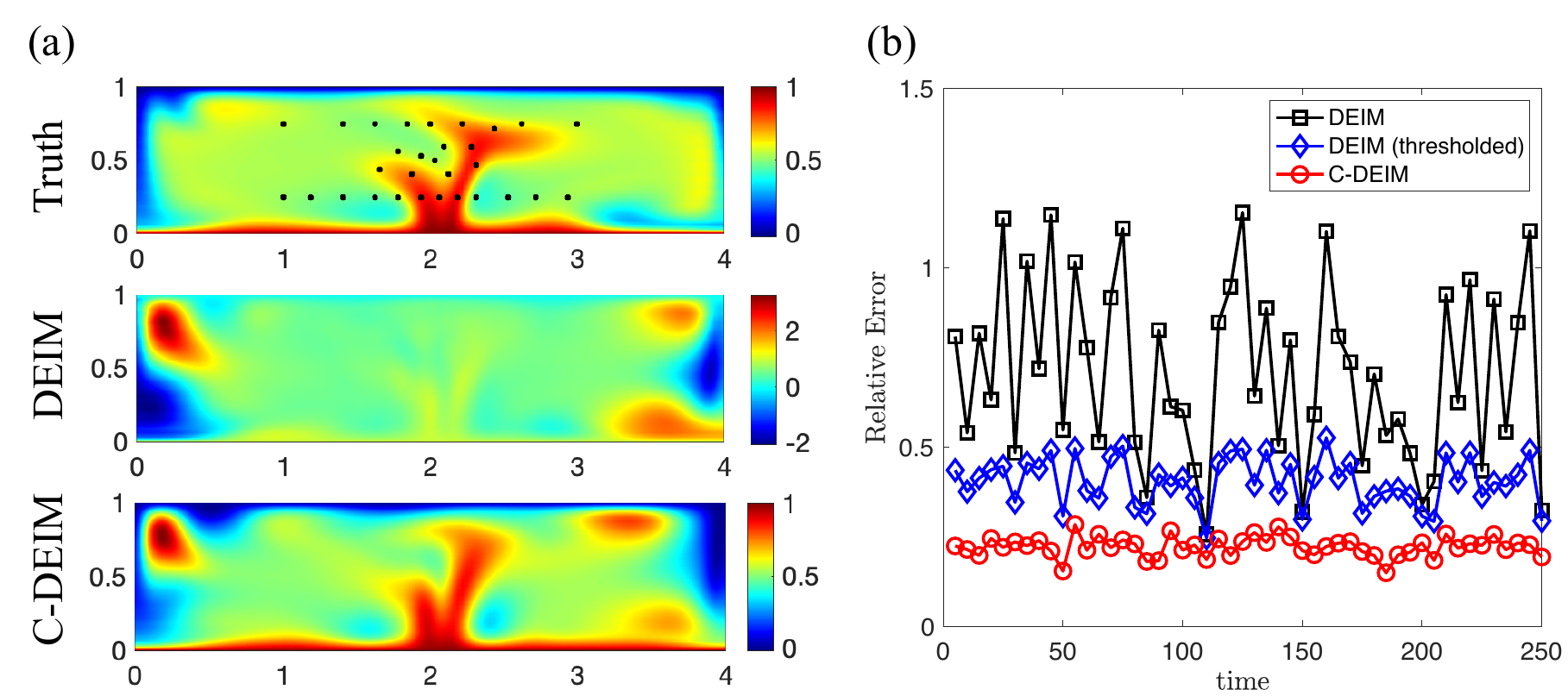}
	\caption{Reconstruction of temperature using DEIM and C-DEIM. (a) A snapshot of temperature $T(\vc x,t)$ from the test data set, and its reconstructions. The black dots mark the location of the sensors. (b) Relative error of the DEIM and C-DEIM reconstructions across all test data.}
	\label{fig:RBC_error}
\end{figure}

We apply DEIM and C-DEIM to the test data to reconstruct the temperature $T(\vc x,t)$ from its partial pointwise sensor measurements.
For these reconstructions, we use $m=30$ POD modes, computed from the training data set, to form the basis matrix $\Phi$.
We use the restricted CPQR method~\eqref{eq:rCPQR} for placing $r=30$ sensors. The sensor locations are restricted to the accessible domain $\mathcal G_a = [1,3]\times [0.25,0.75]$. The resulting sensor locations are marked with black circles in figure~\ref{fig:RBC_error}(a). 
Throughout this section, we use the following input variables for the C-DEIM algorithm: initial penalty parameter $\lambda_{\text{init}} = 10^{-6}$, constraint tolerance $\delta = 10^{-6}$, Newton tolerance $\tau = 10^{-6}$, growth factor $\gamma=10$, and bisection tolerance $\tau_\lambda = 0.1$. The results are not sensitive to small perturbations to these parameters.

We note that the fluid temperature is always bounded between the wall temperatures, i.e., $0\leq T(\vc x,t)\leq 1$. This follows from the boundary conditions~\eqref{eq:RBC_bc} and the maximum principle~\cite{Barrenechea2024}. As a result, a physically consistent reconstruction of the temperature must lie within the interval $[0,1]$. However, as exemplified in figure~\ref{fig:RBC_error}(a), the DEIM reconstruction routinely violates these bounds, reaching values in the range $T>2$ and $T<-1$ which are physically impossible. These violations typically occur in the extrapolation region $\mathcal G \backslash\mathcal G_a$, near the corners of the spatial domain.

In contrast, the C-DEIM reconstructions always remain within the physical bounds $[0,1]$. It is noteworthy that the C-DEIM reconstruction is also inaccurate in the extrapolation domain $\mathcal G \backslash\mathcal G_a$, but nonetheless it respects the physical constraints for the temperature. For this snapshot, the optimal parameter value is $\lambda_{\text{opt}} = 249.08$. Figure~\ref{fig:RBC_error}(b) shows the relative reconstruction error over the entire test data set. The DEIM reconstruction errors are significantly larger than the C-DEIM reconstruction errors. More specifically, the mean relative error for DEIM is around $71\%$, whereas the mean error for C-DEIM is around $22\%$.

Figure~\ref{fig:RBC_error}(b) also shows the relative error for the thresholded DEIM. These reconstructions are obtained by taking the DEIM output and replacing the out-of-bounds values $T>1$ and $T<0$ with $T=1$ and $T=0$, respectively. Even after thresholding, the DEIM reconstructions are worse than the C-DEIM reconstructions. The mean relative error of thresholded DEIM is around $40\%$; almost twice that of C-DEIM.

\subsection{Wildfire Propagation}\label{subsec:fire}

\begin{figure} 
	\centering
	\includegraphics[width=0.4\textwidth]{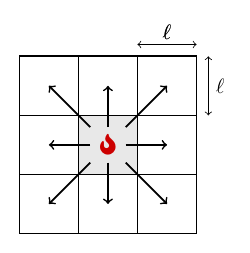}
	\caption{Schematic of wildfire cellular automata model. Fire spreads from a burning cell to neighboring cells in eight directions based on wind-influenced transition rules. Each cell has a side length of $\ell$.}
	\label{fig:ca_model}
\end{figure}

In this section, we apply C-DEIM to reconstruct and forecast the propagation of wildfires. As the base model, we use a cellular automata (CA) approach~\cite{alexandridis2008}, which is detailed in \ref{app:ca}. The spatial domain is discretized into a grid of $N$ square cells, each with side length $\ell$. Each cell $(i,j)$ is assigned a discrete state from the set: \textit{unburned}, \textit{burning}, or \textit{burned}. The fire evolves in discrete time steps $\Delta t$ according to local rules that depend on wind speed and direction, allowing the fire to propagate in eight directions within each cell (see Figure~\ref{fig:ca_model}). These directions are assigned a discrete angle $\theta\in \{0, \pi/4,\pi/2,\cdots,7\pi/4 \}$.

The rate of fire spread from cell $(i,j)$ in direction $\theta$, denoted $R_{ij}(\theta)$, depends on the wind conditions, creating a locally elliptical spread~\cite{peterson2009}. At each time step $\Delta t$, the fire travels a distance $R_{ij}(\theta)\Delta t$ in each direction $\theta$. The cumulative spread in direction $\theta$ after $n$ time steps denoted by $d_{ij}(\theta, n\Delta t)$. When the fire has traveled the distance from the center of the cell to the center of a neighboring cell, the neighboring cell becomes ignited, and the fire progresses in the neighboring cell. 

Each cell $(i,j)$ is assigned a state variable $s_{ij}(t)$, representing the proportion of time the cell has been in the \textit{burning} state up to time $t$. Specifically, 
\begin{equation}
s_{ij}(t)=\begin{cases}
0, \quad \text{if the cell is \textit{unburned} }\\
\frac{t-t^I_{ij}}{t}, \quad \text{if the cell ignited at time $t_{ij}^I\leq t$}.
\end{cases}
\label{eq:fire_cont}
\end{equation}
This yields a bounded state variable $s_{ij}(t)\in[0,1]$ which is a proxy for the portion of the cell that has burned up to time $t$. We form the full state vector $\vc u(t) \in \mathbb{R}^N$ for DEIM by vectorizing $s_{ij}(t)$.

\begin{figure} 
	\centering
	\includegraphics[width=0.5\textwidth]{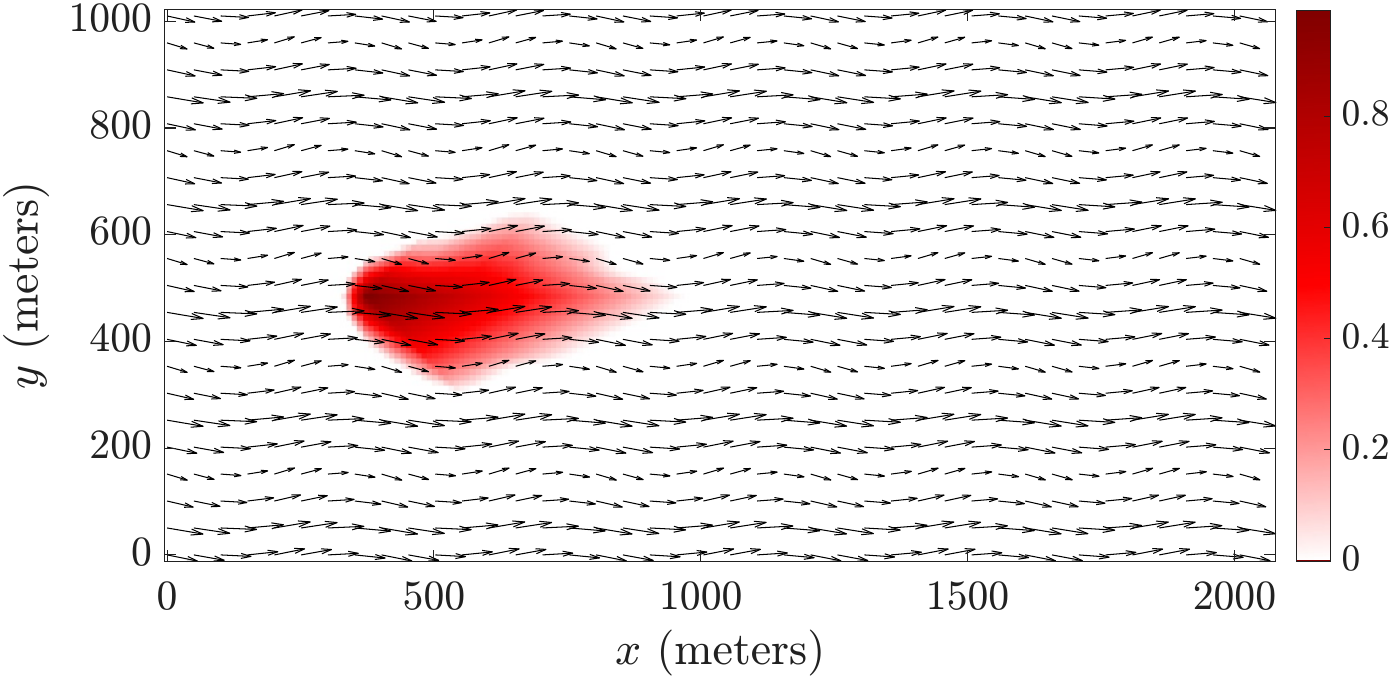}
	\caption{Example simulation of wildfire spread one hour after ignition, using randomly generated wind field  $\vc v(x,y)$~\eqref{eq:wind_v}. The color intensity reflects the proportion of time each cell has been burning.}
	\label{fig:fire_model}
\end{figure}

We simulate variability in fire dynamics by incorporating randomized wind conditions derived from a stream function $\psi(x,y)$ defined by
\begin{equation}
\psi(x,y)=v_0y+\epsilon\left[A\frac{L_x}{\nu}\cos\left(\frac{2\pi \nu x}{L_x}+\varphi_1\right)+B\frac{L_y}{\nu}\sin\left(\frac{2\pi \nu y}{L_y}+\varphi_2\right)\right],
\label{eq:stream}
\end{equation}
where the random amplitudes $A$ and $B$ are drawn from the standard normal distribution $\mathcal{N}(0,1)$ and the random phases $\varphi_1$ and $\varphi_2$ are drawn from a uniform distribution $\mathcal{U}[0,2\pi]$. The amplitude of the perturbation is $\epsilon=0.1$, the spatial frequency is $\nu=5$, and the base wind speed is $v_0=2.5$ m/s. Here, $L_x$ and $L_y$ denote the horizontal and vertical dimensions of the spatial domain, respectively. The wind velocity field $\vc v(x,y)$ at each point $(x,y)\in[0,L_x]\times[0,L_y]$ is given by
\begin{equation}
\vc v(x,y)=\left( \frac{\partial\psi}{\partial y},-\frac{\partial\psi}{\partial x} \right),
\label{eq:wind_v}
\end{equation}
which ensures incompressibility.
Figure~\ref{fig:fire_model} depicts a sample simulation with its corresponding wind velocity field.

We generate 1000 wildfire spread simulations, each with the fire ignition at $(x_{\text{init}},y_{\text{init}})=(380,490)$ meters and a different wind realization.
We record the cell state values $s_{ij}$ one hour after ignition for each simulation. We form the training set from the first 800 of these snapshots, which we use to construct the basis matrix $\Phi$ via POD. The remaining 200 simulations constitute the test data. For each test state (a vector $\vc u$ of $s_{ij}$ values at one hour), we extract sparse sensor measurements $\vc y$ at $r$ sensor locations. For C-DEIM, we input the following parameters: an initial penalty parameter $\lambda_{\text{init}}=10^{-6}$, Newton tolerance $\tau=10^{-10}$, growth factor $\gamma=10$, and constraint tolerance $\delta=10^{-7}$. 

We evaluate C-DEIM under two sensor placement scenarios:
	\begin{enumerate}
		\item Restricted CPQR placement: Sensors are selected with CPQR, but the selection is restricted to a limited, accessible subdomain, mimicking constraints due to cost or physical barriers. Specifically, sensors are confined to three horizontal lines through the center of the domain, located at $y=400$, $500$, and $600$ meters.		
		\item Random placement: Sensors are randomly selected from the set of cells that are actively burning at the time of the measurement, simulating presence-only wildfire data by capturing the spatial occurrence of the fire without explicitly recording its absence. 
			\end{enumerate}

\begin{figure} 
	\centering
	\includegraphics[width=\textwidth]{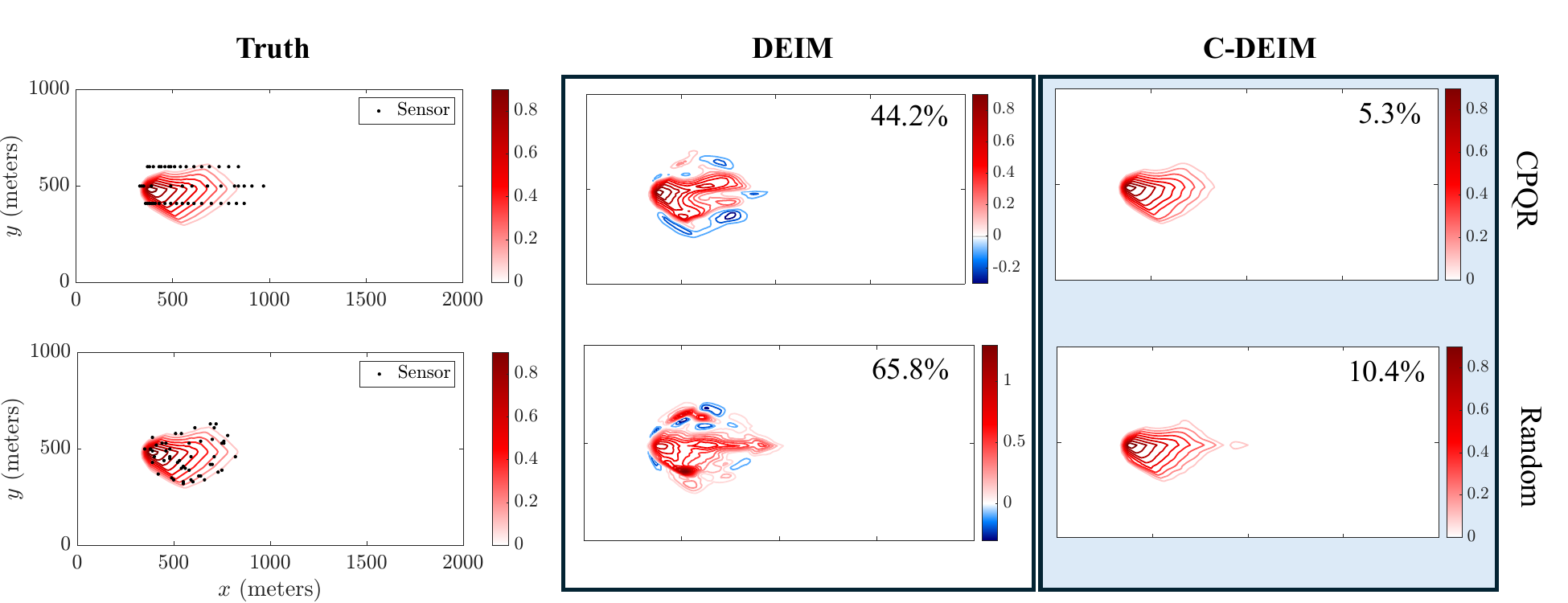}
	\caption{Comparison of wildfire reconstructions at one hour post-ignition using $r=50$ sensors. The relative reconstruction error (expressed as a percentage) is shown in the top right corner of each reconstruction. Restricted CPQR sensor placement: C-DEIM produces a physically valid reconstruction, while DEIM generates negative values (top row). Random sensor placement: Similar trends are observed, with C-DEIM preserving the bounded nature of the state with a significantly smaller relative reconstruction error (bottom row).}
	\label{fig:recon_samples}
	\end{figure}

Figure~\ref{fig:recon_samples} shows representative reconstructions for both sensor placement scenarios using $r=50$ sensors. The DEIM reconstruction includes negative values and values greater than 1, which are not meaningful within the CA model. In contrast, C-DEIM yields a physically reasonable reconstruction, respecting the $[0,1]$ range of the state variable $s_{ij}$. For this function, the optimal penalty parameter is $\lambda_{\text{opt}}=606.25$ for restricted sensor placement and $\lambda_{\text{opt}}=289.84$ for random sensor placement. The C-DEIM reconstruction created with $r=50$ sensors and a growth factor of $\gamma=10$ takes on average approximately 3 minutes more to compute, compared to DEIM. To decrease the computational time, one can increase the growth factor. For example, using $\gamma=100$ reduces the additional computational time of C-DEIM to less than 2 minutes. However, an excessively large growth factor is not advisable either, as it increases the number of bisection iterations required to determine the optimal $\lambda$ (see Section \ref{subsec:update}), potentially increasing the computational time again.

\begin{figure} 
	\centering
	\includegraphics[width=\textwidth]{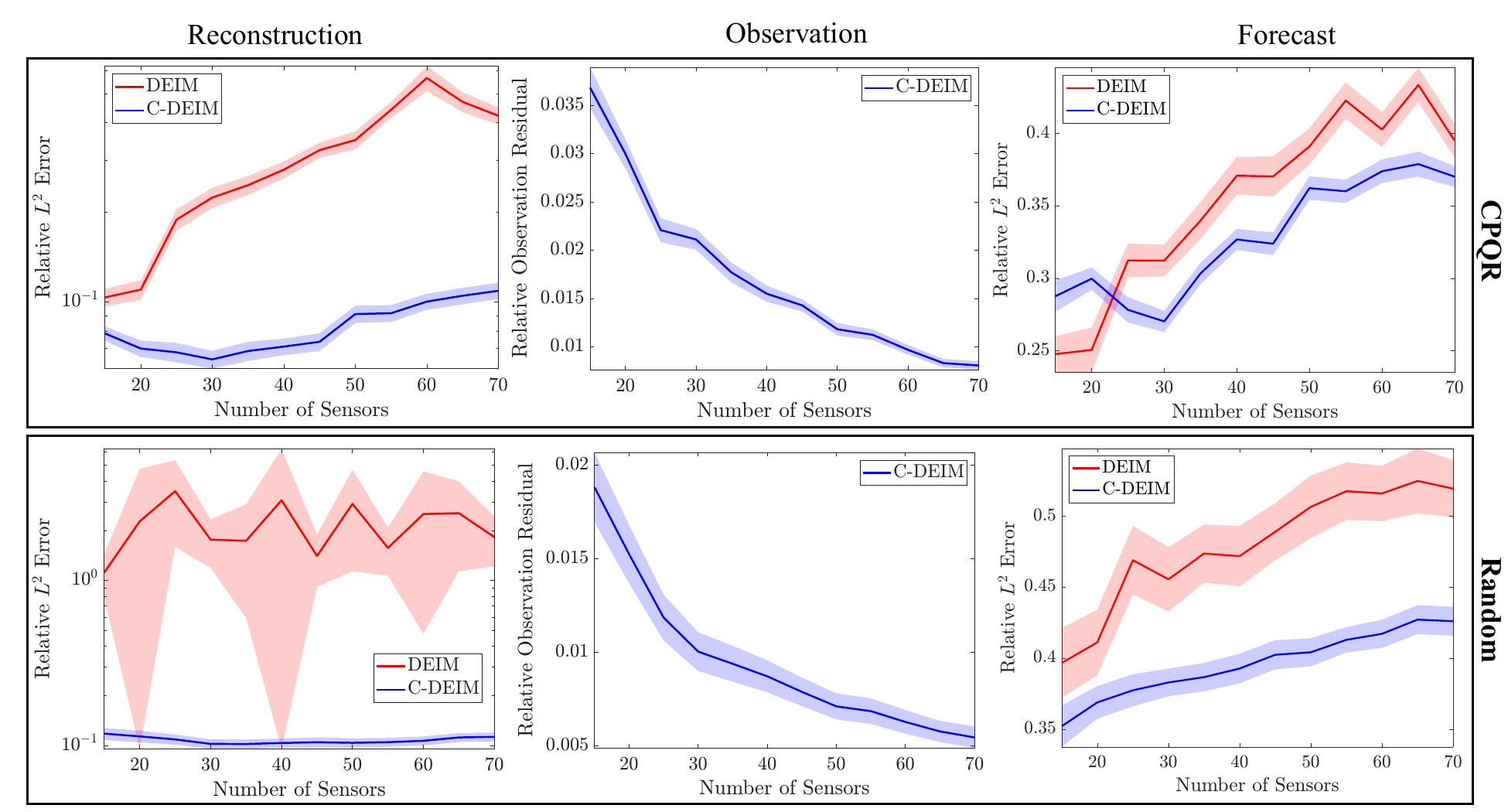}
	\caption{Performance metrics for C-DEIM and DEIM reconstructions and forecasts with increasing number of sensors, restricted CPQR sensor placement (top row), random sensor placement (bottom row). Mean relative reconstruction errors (left), C-DEIM outperforms DEIM in both cases. Mean relative observation residuals (middle). Error decreases with more sensors and remains below 5\%. Mean relative $L^2$ forecast errors, two hours post-ignition (right), using thresholded reconstructions as initial conditions. C-DEIM based forecasts show better accuracy as the number of sensors increases. The shaded regions depict 95\% confidence intervals.}
	\label{fig:fire_errors}
	\end{figure}

Figure~\ref{fig:fire_errors} presents the mean relative reconstruction errors for both methods (left panels) and the mean relative observation residuals (middle panels) across 200 test cases for an increasing number of sensors. Under both restricted CPQR and random sensor placement, C-DEIM achieves a significantly lower mean reconstruction error compared to DEIM. C-DEIM's mean relative observation residual remains low and decreases as more sensors are added, indicating the reconstruction corresponds well with the measurements, despite enforcing constraints. For $r=70$ sensors placed with restricted CPQR, the mean relative reconstruction error for DEIM is 42\%, while C-DEIM achieves a much smaller error of 11\%. For $r=70$ randomly placed sensors, the mean relative reconstruction error for DEIM is 183\% while it is an order of magnitude lower at 11.5\% for C-DEIM. The C-DEIM errors remain about the same regardless of sensor placement method, while the DEIM error is significantly larger when using random sensor placement. This indicates that C-DEIM is not as sensitive to the choice of sensor locations. 

Beyond reconstruction, we assess the efficacy of C-DEIM for improving forecast accuracy. The state reconstruction at one hour post-ignition (obtained using C-DEIM or DEIM) is used as the initial condition for the CA model to predict the fire state at two hours after ignition. Any values in the reconstructed state outside the $[0,1]$ range are thresholded to 0 or 1. Note that the C-DEIM reconstructions require little or no thresholding due to the enforced constraints, whereas the DEIM reconstructions typically require significant thresholding. As the cellular automata progresses with a discrete time step $\Delta t$, we assume an ignition cannot occur in a time interval less than $\Delta t$. Therefore, if $s_{ij}(t)=\frac{t-t^I_{ij}}{t}< \frac{\Delta t}{t}$, we set $s_{ij}(t)=0$ before using the reconstruction as an initial condition. 
 
\begin{figure} 
	\centering
	\includegraphics[width=\textwidth]{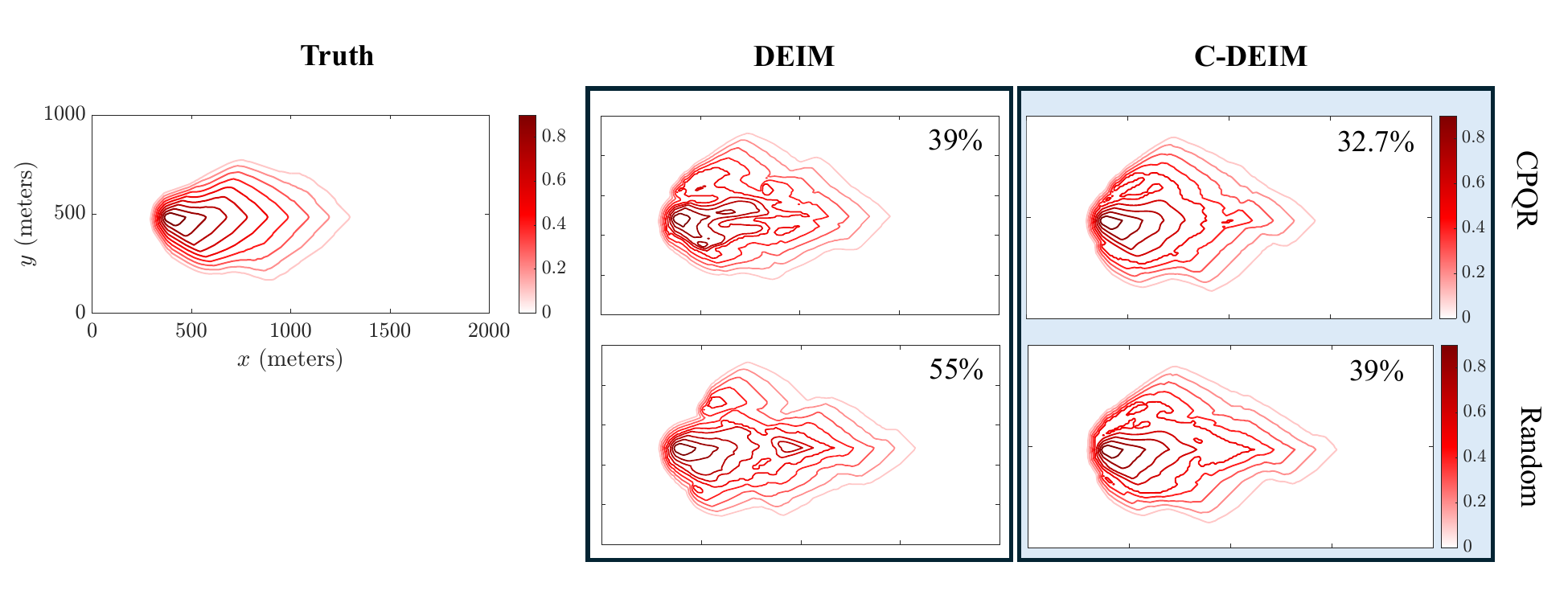}
	\caption{Forecast comparison two hours after fire ignition with $r=50$ sensors, with restricted CPQR sensor placement (top row) and random sensor placement (bottom row). Relative forecast error (expressed as a percentage) is shown in the top right corner of each forecast. C-DEIM-based forecasts better match the true state.}
	\label{fig:forecast_samples}
\end{figure}

Forecast accuracy is evaluated by comparing the predicted state at two hours (based on reconstructions at one hour) to the true state from the full simulation. Figure~\ref{fig:fire_errors} (right panels) shows the mean relative $L^2$ forecast errors. While the DEIM forecasts using the CPQR sensor placement exhibit lower mean relative errors for $r=15$ and $r=20$ sensors, C-DEIM outperforms DEIM in every other case. When sensors are placed randomly, C-DEIM always performs better. A visual comparison of true and forecasted states (figure~\ref{fig:forecast_samples}) further illustrates that C-DEIM-based forecasts provide a better match to the actual fire spread. For restricted CPQR sensor placement, the forecast created with the DEIM reconstruction has a 39\% relative error, while the forecast created with the C-DEIM reconstruction has a smaller error of 32.7\%. For random sensor placement, the C-DEIM forecast has an error of 39\%. The DEIM forecast returns an error of 55\%, which is substantially higher, due to a more elongated and wider initial reconstruction (figure \ref{fig:recon_samples}, bottom center panel). 

In summary, the forecasts corresponding to restricted CPQR sensor placement are more accurate than those from random sensor placement. However, in practice, determination and deployment of sensors at optimal locations is often not possible. As such, random sensor placement mimics the operational constraints in data gathering. In this case, the forecasts based on C-DEIM reconstructions are more accurate than their DEIM-based counterparts (cf. figure~\ref{fig:fire_errors}, right panels).

\section{Conclusions}\label{sec:concl}
DEIM is a simple, yet powerful, method for reconstructing a field from its partial pointwise observations. However, DEIM reconstructions can violate prescribed upper and lower bounds of the quantity of interest. This may limit the utility of the reconstruction in downstream tasks, such as forecasting or control.

Here, we developed C-DEIM to address this issue. By employing soft penalty constraints, C-DEIM explicitly incorporates the prescribed bounds into the reconstruction process. We proved that the solution of the penalized minimization problem converges to the constraint set, as the penalty parameter increases (Theorem~\ref{thm:convergence_lam}), providing a theoretical guarantee that the C-DEIM reconstruction will respect the prescribed bounds.

We also derived a theoretical upper bound on the observation residual of the C-DEIM reconstruction (Theorem~\ref{thm:obs_err}). This upper bound informs the selection of the optimal penalty parameter, using a bisection method. As a result, the constraints can be satisfied to arbitrary accuracy, while simultaneously ensuring that the observation residual remains as small as possible. Furthermore, we designed C-DEIM based on the above theoretical results and proved that it converges in a finite number of iterations.

We demonstrated the efficacy of C-DEIM through three numerical examples:  random harmonic functions, the Rayleigh–B\'enard convection, and a wildfire propagation model. The results showed that C-DEIM produces significantly improved reconstruction accuracy compared to DEIM, particularly when sensor placement is restricted to a subset of the spatial domain or when the sensors are distributed randomly.

Several avenues for future research exist:
\begin{itemize}
	\item Investigating optimal sensor placement strategies tailored for C-DEIM would be valuable.
	\item The theoretical error bound~\eqref{eq:obs_ub} on the observation residual provides a useful analytical result, but numerical results suggest it may be pessimistic, and a tighter error bound may be attainable. Furthermore, extending this result (Theorem~\ref{thm:obs_err}) to the overdetermined regime $r>m$ is desirable.
	\item Our current approach reconstructs the quantity of interest from a single snapshot of observations. For dynamical systems, incorporating sequential time series data enhances reconstruction and forecasting accuracy. This sequential approach was recently developed for DEIM~\cite{farazmand2024,Farazmand2025}; its development in the context of C-DEIM will be valuable.
\end{itemize}  

Finally, C-DEIM can be readily applied to experimental data. Its demonstrated robustness under restricted and random sensor placement suggests promising potential for practical use in applications. A manuscript by the authors on the application of C-DEIM to the 2023 Maui wildfire is forthcoming.

\section*{Data Availability}
No data was used for the research described in this article.

\section*{Funding}
This work was supported by the National Science Foundation RTG Program, through the award DMS-2342344, and the Algorithms for Threat Detection (ATD) program, through the award DMS-2220548.

\appendix
\section{Proof of Theorem \ref{thm:convergence_lam}}
\label{app:conv}
We first prove a number of results regarding Lyapunov-like functions which are subsequently used to prove Theorem~\ref{thm:convergence_lam}. Consider a differential equation,
\begin{equation}
	\dot{\vc x} = \vc f(\vc x),
	\label{eq:ode}
\end{equation}
where the vector field $\vc f:\mathbb R^n\to\mathbb R^n$ is Lipschitz continuous. We denote the solution map corresponding to~\eqref{eq:ode} by $\fm^t:\mathbb R^n\to\mathbb R^n$ so that $\fm^t(\vc x_0)$ denotes the trajectory of the system emanating from the initial condition $\vc x_0$.
A Lyapunov function is a map $V:\mathbb R^n\to\mathbb R$ such that $V(\vc x)>0$ for all $\vc x$ except at a fixed point $\vc x_\ast$ where $V(\vc x_\ast) =0$.
Furthermore, $V$ is non-increasing along the trajectories of~\eqref{eq:ode}, so that $\dot V(\vc x)\leq 0$. Lyapunov functions are routinely used to prove stability of fixed points of dynamical systems. Here, we extend classical results to Lyapunov-like functions which vanish over a simply connected set instead of a single fixed point. 

\begin{defn}
	We denote the interior of the level sets of $V$ by $S_v$. More precisely, for any $v>0$, we define
	$$S_v= \{\vc x\in\mathbb R^n: V(\vc x)<v\}.$$
	For $v=0$, we define $S_0$ to be the corresponding level set,
	$$ S_0 = \{\vc x\in\mathbb R^n: V(\vc x) =0\}.$$
\end{defn}

\begin{lem}\label{lem:Lyap1}
	Assume that there exists a Lyapunov-like function $V:\mathbb R^n\to\mathbb R$ for system~\eqref{eq:ode} which satisfies the following conditions:
	\begin{enumerate}
		\item $V(\vc x) = 0 $ for all $\vc x\in S_0$.
		\item $V(\vc x) > 0$ for all $\vc x\in\mathbb R^n\backslash S_0$. 
		\item $\dot V(\vc x)\leq 0$ for all $\vc x\in\mathbb R^n$.
	\end{enumerate}
	
	Then the sets $S_v$ are positively invariant for all $v\leq 0$, i.e., if $\vc x_0\in S_v$ then $\fm^t(\vc x_0)\in S_v$ for all $t\geq 0$.
\end{lem}
\begin{proof}
	Let's assume to the contrary that there exists $\vc x_0\in S_v$ such that $\fm^T(\vc x_0)\notin S_v$ for some $T>0$.
	Since the trajectories are continuous and $V$ is continuous, there must exist $0<t<T$ such that $\fm^t(\vc x_0) \in \partial S_v$, where $\partial S_v$ denotes the boundary of $S_v$. This implies $V(\vc x_0)<v<V(\fm^t(\vc x_0))= v$. This contradicts $\dot V\leq 0$.
\end{proof}

\begin{lem}\label{lem:Lyap2}
	Assume that there exists a Lyapunov-like function $V:\mathbb R^n\to\mathbb R$ for system~\eqref{eq:ode} which satisfies the following conditions:
	\begin{enumerate}
		\item $V(\vc x) = 0 $ for all $\vc x\in S_0$.
		\item $V(\vc x) > 0$ for all $\vc x\in\mathbb R^n\backslash S_0$. 
		\item $\dot V(\vc x)< 0$ for all $\vc x\in\mathbb R^n\backslash S_0$.
	\end{enumerate}
	
	Then, for all $\vc x_0\in\mathbb R^n$, we have
	\begin{equation}
		\lim_{t\to\infty}\id(\fm^t(\vc x_0), S_0) = 0.
	\end{equation}
\end{lem}

\begin{proof}
\textbf{Part 1.} We first show that for any $\vc x_0$ and any $v>0$, there exist $T>0$ such that $\fm^t(\vc x_0)\in S_v$ for all $t\geq T$.
Note that Lemma~\ref{lem:Lyap1} implies that if the trajectory $\fm^t(\vc x_0)$ enters $S_v$, it stays in $S_v$ for all times.
More precisely, if $\fm^T(\vc x_0)\in S_v$ for some $T\geq 0$ then $\fm^t(\vc x_0)\in S_v$ for all $t\geq T$.

Now assume to the contrary that there exist $\vc x_0$ and $v>0$, such that $\fm^t(\vc x_0)\notin S_v$ for all $t\geq 0$.
Let $V(\vc x_0) = v'>v$. Then $\fm^t(\vc x_0)\in \overline{S_{v'}}\backslash S_v$ for all $t\geq 0$. Since the set  $\overline{S_{v'}}\backslash S_v$ is compact, the function $\dot V$ attains its minimum somewhere on this set. Furthermore, since $\dot V<0$ on this set, the minimum is negative. Therefore, there exists $\mu>0$ such that 
$$\dot V(\fm^t(\vc x_0))\leq -\mu,\quad \forall t\geq 0.$$
Integrating in time, we obtain
$$ V(\fm^t(\vc x_0))\leq V(\vc x_0)-\mu t.$$
This implies that there exists $t>0$ large enough so that $V(\fm^t(\vc x_0))<0$. This contradicts $V$ being non-negative.

\textbf{Part 2.} Consider a monotonically decreasing sequence $v_1>v_2>v_3>\cdots>0$ such that $\lim_{n\to\infty} v_n = 0$. The sets $S_{v_n}$ are nested $S_{v_1}\supset S_{v_2}\supset\cdots\supset S_0$. Furthermore, $S_{v_n}$ converges to $S_0$ as $n$ tends to infinity, i.e., 
\begin{equation}
	\lim_{n\to\infty} \sup_{\vc x\in S_{v_n}} \id (\vc x,S_0)=0.
\end{equation}
In particular, this implies that there exists a decreasing sequence $\epsilon_n$ such that $\lim_{n\to\infty} \epsilon_n= 0$ and $\sup_{\vc x\in S_{v_n}} \id (\vc x,S_0)<\epsilon_n$.

From Part 1, we conclude that for any $\vc x_0$, there exists $T_n$ such that 
$\fm^t(\vc x_0)\in S_{v_n}$ for all $t\geq T_n$. This implies 
\begin{equation}
	0\leq \liminf_{t\to\infty} \id (\fm^t(\vc x_0),S_0)\leq \limsup_{t\to\infty} \id (\fm^t(\vc x_0),S_0)\leq \epsilon_n.
\end{equation}
Since $\lim_{n\to\infty}\epsilon_n=0$, we obtain the desired result that $\lim_{t\to\infty} \id (\fm^t(\vc x_0),S_0)=0.$
\end{proof}

We define the total penalty $\mathcal{P}(\balpha)$ as the Lyapunov-like function and the corresponding level set $\mathcal{S}$  where $\mathcal{P}(\balpha)=0$ for all $\balpha \in \mathcal{S}$. The penalty parameter $\lambda$ plays the role of time $t$. Then, the following theorem guarantees that the shortest distance between the minimizer $\balpha(\lambda)$~\eqref{eq:alpha_lam} and constraint set $\mathcal{S}$ goes to zero, as the penalty parameter $\lambda$ increases to infinity.\\

\begin{proof}[Proof of Theorem \ref{thm:convergence_lam}]
The solution $\balpha(\lambda)$ for a given $\lambda\geq0$ minimizes the penalized cost function \eqref{eq:obj_f}. The first-order optimality condition is given by
	\begin{equation}
	\nabla\err_o(\balpha(\lambda))+\lambda\nabla\mathcal{P}(\balpha(\lambda))=\mathbf{0}.
	\label{eq:alpha_deriv}
	\end{equation}
	For the trivial case, where the solution for $\lambda=0$ is within the constraint set ($\balpha(0)\in \mathcal{S}$), the total penalty $\mathcal{P}(\balpha(0))=0$. Therefore, $\nabla \mathcal{P}(\balpha(\lambda))$ vanishes, so $\balpha(\lambda)\in\mathcal{S}$ for all $\lambda>0$. 
	
	For $\balpha(0)\not\in\mathcal{S}$, we analyze the behavior of the solution $\balpha({\lambda})$ as $\lambda \to \infty$. Differentiating the optimality condition \eqref{eq:alpha_deriv} with respect to $\lambda$, using the chain rule for terms involving $\balpha(\lambda)$, yields,
	\begin{equation}
	\nabla^2\err_o(\balpha(\lambda))\balpha'(\lambda)+\lambda\nabla^2\mathcal{P}(\balpha(\lambda))\balpha'(\lambda)+\nabla\mathcal{P}(\balpha(\lambda))=\mathbf{0}.
	\label{eq:opt_cond}
	\end{equation}
	where $\balpha'(\lambda)=\frac{\text{d}\balpha}{\text{d}\lambda}$. The Hessian $\nabla^2\err_o(\balpha)$ is positive semidefinite for all $\balpha \in \mathbb{R}^m$. The diagonal matrix $D_{ii}(\Phi\balpha)=p''([\Phi\balpha]_i)$ has nonnegative entries, since $p$ is convex. For $\balpha\not\in\mathcal{S}$, $p''([\Phi\balpha]_i)>0$ for at least one coordinate $i$. The matrix $\Phi$ has orthonormal columns, so for $\lambda>0$, the Hessian of the penalty term, $\nabla^2\mathcal{P}(\balpha(\lambda))=\Phi^\top D(\Phi\balpha(\lambda))\Phi$ is positive definite for all $\balpha\not\in\mathcal{S}$. Thus, for $\lambda>0$, the matrix $[\nabla^2\err_o(\balpha(\lambda))+\lambda\nabla^2\mathcal{P}(\balpha(\lambda))]$ is the sum of a positive semidefinite matrix and a positive definite matrix, which is positive definite and therefore nonsingular. Thus, we can solve \eqref{eq:opt_cond} for $\balpha'(\lambda)$,
	\begin{equation}
	\balpha'(\lambda)=-[\nabla^2\err_o(\balpha(\lambda))+\lambda\nabla^2\mathcal{P}(\balpha(\lambda))]^{-1}\nabla \mathcal{P}(\balpha(\lambda)).
	\end{equation}
	Now, consider the derivative of the total penalty $\mathcal{P}(\balpha(\lambda))$ with respect to $\lambda$ using the chain rule, 	
	\begin{equation}
	\frac{\text{d}}{\text{d}\lambda}\mathcal{P}(\balpha(\lambda))=-\nabla\mathcal{P}(\balpha(\lambda))^\top[\nabla^2\err_o(\balpha(\lambda))+\lambda\nabla^2\mathcal{P}(\balpha(\lambda))]^{-1}\nabla\mathcal{P}(\balpha(\lambda)).
	\label{eq:deriv_lam}
	\end{equation}
	Since the matrix $[\nabla^2\err_o(\balpha(\lambda))+\lambda\nabla^2\mathcal{P}(\balpha(\lambda))]$ is positive definite, its inverse is also positive definite. Therefore, if  $\nabla\mathcal{P}(\balpha(\lambda)) \ne \mathbf{0}$, $\frac{\text{d}}{\text{d}\lambda}\mathcal{P}(\balpha(\lambda)) < 0$. The total penalty $\mathcal{P}(\balpha(\lambda))$ strictly decreases as $\lambda$ increases, as long as $\nabla\mathcal{P}(\balpha(\lambda)) \ne \mathbf{0}$. The derivative \eqref{eq:deriv_lam} is zero only when $\nabla\mathcal{P}(\balpha(\lambda)) = \mathbf{0}$. However, $\nabla\mathcal{P}(\balpha(\lambda)) = \mathbf{0}$ if and only if $\balpha \in \mathcal{S}$.
	
Then, we can use $\mathcal{P}(\balpha(\lambda))$ as a Lyapunov-like function that satisfies three conditions: $\mathcal{P}(\balpha(\lambda))=0$ for all $\balpha \in \mathcal{S}$, $\mathcal{P}(\balpha(\lambda))>0$ for all $\balpha\in \mathbb R^m\backslash \mathcal{S}$, and $\frac{\text{d}}{\text{d}\lambda}\mathcal{P}(\balpha(\lambda))<0$ for all $\balpha \in\mathbb R^m\backslash \mathcal{S}$. According to Lemma~\ref{lem:Lyap2}, 
\begin{equation}
\lim_{\lambda\to \infty}\text{d}(\balpha(\lambda),\mathcal{S})=0, 
\end{equation}
which implies that as $\lambda \to \infty$, $\balpha(\lambda)$ converges to the constraint set $\mathcal{S}$.
\end{proof}

\section{Proof of Theorem~\ref{thm:conv}}
\label{app:newton}
We first state a lemma regarding the Lipschitz continuity of diagonal matrices, formed from a Lipschitz function applied element-wise to a vector. 

\begin{lem}\label{lem:lip_cont}
	Let $\mathbf{x}\in\mathbb{R}^{n}$ with components $x_i$, and let $D(\vc x)\in \mathbb{R}^{n\times n}$ be a diagonal matrix such that $D_{ii}(\vc x)=h(x_i)$ for some function $h : \mathbb{R} \to \mathbb{R}$. If $h(x)$ is Lipschitz continuous, then $D(\vc x)$ is also Lipschitz continuous. \\
	\end{lem}
	\begin{proof}
	If $h(x)$ is Lipschitz continuous, then there exists a constant $L$, such that
	\begin{equation*}
	|h(x)-h(y)|\leq L |x-y|,
	\end{equation*}
	for all $x,y\in \mathbb{R}$.
	The spectral norm of a diagonal matrix is the maximum absolute value of the diagonal entries. Thus, for any $\vc x, \vc y \in \mathbb{R}^{n}$,
	\begin{equation*}
	\|D(\vc x)-D(\vc y))\|_2=\max_{i}|h(x_i)-h(y_i)|.
	\end{equation*}
	Let $j$ be an index such that $\max_{i}|h(x_i)-h(y_i)|=|h(x_j)-h(y_j)|$.
	Then, using the Lipschitz continuity of $h$,
	\begin{equation*}
	\|D(\vc x)-D(\vc y)\|_2=|h(x_j)-h(y_j)| \leq L |x_j-y_j|\leq L |\vc x-\vc y|.
	\end{equation*}
	Thus, $D(\mathbf{x})$ is Lipschitz continuous with constant $L$.
	\end{proof}

We use Lemma~\ref{lem:lip_cont} to show that the Hessian of $f_\lambda(\balpha)$ is positive definite and Lipschitz continuous, satisfying conditions for Newton convergence. 

\begin{proof}[Proof of Theorem~\ref{thm:conv}]
For $\lambda >0$, the penalized cost function $f_\lambda(\balpha)$ is nonnegative and convex, ensuring the existence of a minimizer where $\nabla f_\lambda(\balpha(\lambda))=\mathbf{0}$. The Newton iterations satisfy the conditions of convergence for the minimizer $\balpha(\lambda)$ if the Hessian $\nabla^2f_\lambda(\balpha)$ is positive definite and Lipschitz continuous in a neighborhood of $\balpha(\lambda$). We will show that the sufficient conditions are satisfied globally under the theorem's assumptions.

The Hessian of $f_\lambda(\balpha)$ is given by
\begin{equation}
\nabla^2 f_\lambda(\balpha)=\Theta^{\top}\Theta+\lambda\Phi^\top D(\Phi\balpha)\Phi,
\label{eq:hessian}
\end{equation}
where the diagonal matrix $D_{ii}(\Phi \balpha) =p''([\Phi \balpha]_i)$ is formed from the second derivatives of the component-wise penalty function $p''(\tilde{u}_i$) evaluated at each component of the reconstruction.  

Note that the matrix $\Theta^\top\Theta$ is positive semidefinite. Furthermore, the matrix $\Phi^TD(\Phi\balpha)\Phi$ is positive definite for all $\balpha\notin\mathcal S$. The latter follows from the fact that the matrix $D(\Phi\balpha)$ is diagonal with positive entries $D_{ii}(\Phi\balpha)=p''([\Phi\balpha]_i)>0$. Positive definiteness of $\Phi^TD(\Phi\balpha)\Phi$ follows from the fact that $\Phi$ has orthonormal columns. For $\lambda>0$, the Hessian \eqref{eq:hessian} is the sum of a positive semidefinite matrix and a positive definite matrix, which implies $\nabla^2 f_\lambda(\balpha)$ is positive definite for all $\balpha\notin\mathcal S$. This also implies the Hessian is nonsingular. 

To prove Lipschitz continuity of the Hessian, we note that
\begin{align}
\|\nabla^2 f_\lambda(\mathbf{x})-\nabla^2 f_\lambda(\mathbf{y})\|_2&=\|\Theta^\top \Theta +\lambda \Phi^\top D(\Phi\mathbf{x})\Phi-\Theta^\top \Theta -\lambda \Phi^\top D(\Phi\mathbf{y})\Phi\|_2 \notag \\ 
&= \lambda \| \Phi^{\top}[D(\Phi\mathbf{x})-D(\Phi\mathbf{y})]\Phi \|_2 \notag \\
&\leq  \lambda \| \Phi^{\top}\|_2 \|D(\Phi\mathbf{x})-D(\Phi\mathbf{y})\|_2 \|\Phi \|_2.
\end{align}
The basis matrix $\Phi$ has orthonormal columns, implying  $\|\Phi \|_2=\|\Phi^\top\|_2=1$. Therefore, we obtain
\begin{equation}
\|\nabla^2 f_\lambda(\mathbf{x})-\nabla^2 f_\lambda(\mathbf{y})\|_2 \leq \lambda \|D(\Phi\mathbf{x})-D(\Phi\mathbf{y})\|_2.
\end{equation}
Since $p''$ is Lipschitz continuous, the diagonal matrix $D(\Phi\vc x)$ is Lipschitz continuous by Lemma~\ref{lem:lip_cont}. Let $L$ be the Lipschitz constant for $p''$. Then, $ \|D(\Phi \vc x)-D(\Phi \vc y)\|_2 \leq L |\Phi\vc x-\Phi \vc y|$. Furthermore, $|\Phi \vc x -\Phi \vc y| =|\Phi(\vc x -\vc y)|=|\vc x-\vc y|$, so
\begin{equation}
 \|\nabla^2 f_\lambda(\mathbf{x})-\nabla^2 f_\lambda(\mathbf{y})\|_2 \leq \lambda L |\mathbf{x}-\mathbf{y}|.
 \end{equation} 
 Thus, $\nabla^2 f_\lambda(\balpha)$ is Lipschitz continuous with constant $\lambda L$. 

Therefore, the Hessian $\nabla^2 f_\lambda(\balpha)$ is nonsingular and Lipschitz continuous for any $\balpha \notin\mathcal S$ and $\lambda>0$, satisfying the sufficient conditions for local convergence of Newton's method.
\end{proof}

\section{Cellular Automata for Wildfire Propagation}
\label{app:ca}
Cellular automata (CA) models simulate fire spread and capture their highly nonlinear dynamics. These models account for environmental factors such as wind speed and direction, topography, and spatially varying spread rates. Clarke et al.~\cite{clarke1994} first applied these models for practical use, while Karafyllidis and Thanailakis~\cite{karafyllidis1997} presented an original model for mathematical analysis. Subsequent models improved accuracy by incorporating more realistic spread mechanisms~\cite{berjak2002,encinas2007,ghisu2015}. Peterson et al.~\cite{peterson2009} integrated an elliptical spread model~\cite{albini1980} into the CA framework, enabling the simulation of wind-elongated fire fronts. Additional modifications include the use of hexagonal cells~\cite{trunfio2004}, neural networks to derive adaptive transition rules~\cite{zheng2017}, and probabilistic outcomes to reflect uncertainty in fire behavior~\cite{trucchia2020}. We choose to build our model based on the elliptical spread CA model~\cite{peterson2009} for its computational simplicity and its quantifiable time step, which is expressed in terms of the spread rate. 

The CA used in Section~\ref{subsec:fire} discretizes the computational domain $(x,y)\in[0,L_x]\times [0,L_y]$ into a grid of $N$ square cells of length $\ell$. We use $L_x=2000$ meters and $L_y=1000$ meters with side length $\ell=10$ meters. Each cell is assigned one of three discrete states: \textit{unburned}, \textit{burning}, or \textit{burned down}; these assignments are updated as the fire evolves. Fire evolution is governed by local transition rules applied at discrete time steps $\Delta t$. These rules depend primarily on the wind conditions, which are derived from a randomized stream function $\psi(x,y)$ defined in \eqref{eq:stream}. The wind velocity field reads
\begin{equation}
\mathbf{v}(x,y)=\begin{bmatrix}
v_0+\epsilon\left[2\pi B \cos(\frac{2\pi \nu y}{L_y}+\varphi_2)\right]\\
\epsilon\left[2\pi A \sin(\frac{2\pi \nu x}{L_x}+\varphi_1)\right]
\end{bmatrix}.
\label{eq:wind}
\end{equation}
The parameters, $A,B, \varphi_1,\varphi_2, \epsilon$, $\nu$ and $v_0$ are detailed in Section~\ref{subsec:fire}. The local wind velocity in cell $(i,j)$ is denoted by $\vc v_{ij}$. 

The fire spread rate is calculated for each cell and direction based on wind conditions. The maximum potential fire spread rate $R^{\max}_{ij}$ in each cell is calculated using the 10\% rule of thumb~\cite{cruz2019} and a small diffusion term ($5\cdot 10^{-3}$) to allow for minimal spread in the absence of wind, 
\begin{equation}
R^{\max}_{ij}=0.1|\vc v_{ij}|+5\cdot 10^{-3}.
\end{equation}
The spread rate $R_{ij}(\theta)$ in each direction $\theta$ within cell $(i,j)$ is then calculated to ensure the fire spreads in an elliptical shape~\cite{albini1980},
\begin{equation}
R_{ij}(\theta)=R^{\max}_{ij}\frac{1-E_{ij}}{1-E_{ij}\cos(\theta-\theta^v_{ij})},
\end{equation}
where $\theta^v_{ij}$ is the direction of the wind, and $E_{ij}$ is the eccentricity of the fire ellipse. This eccentricity is calculated through the length-to-width ratio $\rho_{ij}$ in cell $(i,j)$,
\begin{equation}
E_{ij}=\sqrt{1-\left( \frac{1}{\rho_{ij}}\right)^2}.
\end{equation}
The length-to-width ratio $\rho_{ij}$ for the ellipse in cell $(i,j)$  is empirically related to wind speed magnitude by $\rho_{ij}=1+0.5592|\vc v_{ij}|$~\cite{rothermel1991}.

Following~\cite{peterson2009}, the maximum simulation time step $\Delta t_{\max}$ is calculated using,
\begin{equation}
\Delta t_{\max}=\frac{\ell}{\max_{\theta,i,j} [R_{ij}(\theta)]},
\end{equation}
to ensure that the fire cannot spread by more than one cell length $\ell$ in any direction during a single time step. To simulate the model for a specified time $T$, we ensure that the time step $\Delta t$ divides evenly into $T$ by using a slightly smaller time step as follows,
\begin{equation}
\Delta t=\frac{T}{\left\lceil T/\Delta t_{\max}\right\rceil}.
\label{eq:time_step}
\end{equation}

Fire propagates by accumulating distance in each direction from a burning cell. For cell $(i,j)$ that has been burning for $n$ time steps, the cumulative spread distance in direction $\theta$ is tracked. When this cumulative distance $d_{ij}(\theta,n\Delta t)$ exceeds the distance from the center of the cell to the center of neighboring cell in direction $\theta$ (which is $\ell$ for direct neighbors and $\ell\sqrt{2}$ for diagonal neighbors), the neighboring cell ignites and transitions to the \textit{burning} state. Any accumulated distance exceeding the cell-to-cell distance is transferred to the newly ignited cell as spillover fire distance to start its own spread accumulation. A cell in the \textit{burning} state transitions to the \textit{burned down} state once all eight of its neighboring cells are \textit{burning}. The fire simulation starts with one initially ignited cell.

At the end of the simulation, the discrete cell states (\textit{unburned}, \textit{burning}, \textit{burned down}) are mapped to a continuous state variable $s_{ij}(t)\in[0,1]$, as defined in \eqref{eq:fire_cont}, representing the proportion of each cell that has burned. The full state vector $\vc u$ for DEIM is formed by vectorizing these $s_{ij}$ values across the grid. 

\subsection{Model Implementation}
The steps for implementing the CA model are detailed below.
\begin{enumerate}
\item Set domain size $[0,L_x]\times[0,L_y]$ and cell side length $\ell$. Set each cell $(i,j)$ to the \textit{unburned} state.
\item Generate the wind velocity field $\vc v(x,y)$ from the randomized stream function $\psi(x,y)$.
\item Compute the length-to-width ratio $\rho_{ij}$, ellipse eccentricity $E_{ij}$, and spread rates $R_{ij}(\theta)$ for each cell $(i,j)$. 
\item Determine the simulation time step $\Delta t$ according to \eqref{eq:time_step}. 
\item Initialize the ignition by setting the cell $(x_{\text{init}},y_{\text{init}})$ to \textit{burning} at time $t=0$.
\item For every $n$-th time step, calculate the spread distance $d_{ij}(\theta, n\Delta t)$ for each direction $\theta$ in each cell $(i,j)$. Evolve the cell states according to the transition rules, tracking the ignition time $t^I_{ij}$ of each \textit{burning} cell. 
\item Convert discrete cell states to continuous values $s_{ij}(t)$ using the current simulation time $t$ and ignition time $t^I_{ij}$.
\end{enumerate}
This framework can be used to simulate the wildfire spread for any simulation time $T$. The parameters and quantities used in the cellular automata model are summarized in Table \ref{tab:ca}.
\begin{table}[ht]
\centering
\begin{tabular}{|c|l|}
\hline
\textbf{Symbol} & \textbf{Description} \\
\hline \hline
$N$ & Number of cells \\ \hline
$[0,L_x]\times [0,L_y]$ & Domain size \\ \hline
$\ell$ & Cell side length \\ \hline
$\Delta t$ & Simulation time step\\ \hline
$T$ & Simulation time \\ \hline
$t^I_{ij}$ & Ignition time for cell $(i,j)$ \\ \hline
$\vc v_{ij}$ & Wind velocity vector in cell $(i,j)$\\ \hline
$\theta^v_{ij}$ & Wind direction in cell $(i,j)$ \\ \hline
$\psi (x,y)$ & Randomized stream function\\ \hline
$A,B$ & Random amplitudes for $\psi$\\ \hline
$\varphi_1,\varphi_2$ & Random phases for $\psi$\\ \hline
$v_0$ & Base wind speed\\ \hline
$\epsilon$ & Amplitude of perturbation for $\psi$\\ \hline
$\nu$ & Spatial frequency for $\psi$\\ \hline
$\rho_{ij}$ & Ellipse length-to-width ratio in cell $(i,j)$  \\ \hline
$E_{ij}$ & Ellipse eccentricity in cell $(i,j)$  \\ \hline
$R_{ij}(\theta)$ & Fire spread rate in cell $(i,j)$ and direction $\theta$  \\ \hline
$d_{ij}(\theta,n\Delta t)$ & Distance of spread from center of cell $(i,j)$ in direction $\theta$ after $n\Delta t$  \\ \hline
$s_{ij}(t)$ & Proportion of time cell $(i,j)$ has been ignited at time $t$ \\
\hline
\end{tabular}
\caption{The parameters and quantities used in the cellular automata model.}
\label{tab:ca}
\end{table}

Wildfires propagate through three primary modes of heat transfer: convection, radiation, and mass transport~\cite{hoffman2021}. We model the convective and radiative heat transfer through the spread rate $R_{ij}(\theta)$. This local spread occurs from the ignition point to the neighboring areas. In contrast, mass transport occurs in the form of spotting, when burning materials are carried downwind. These burning materials can ignite new fires far from the main fire front, called spot fires. Some wildfire cellular automata models incorporate stochastic spotting processes~\cite{alexandridis2008}. In this work, however, we focus exclusively on the local spread. Incorporating spotting into the reconstruction and predictive framework is a promising direction for future research.

\end{document}